\documentclass[Afour,sageh,times]{sagej}

\usepackage{moreverb,url}

\usepackage[colorlinks,bookmarksopen,bookmarksnumbered,citecolor=red,urlcolor=red]{hyperref}

\usepackage[normalem]{ulem}

\usepackage{tikz}
\usepackage{pgfplots}
\usepackage{caption}
\usepackage{subcaption}

\newcommand\BibTeX{{\rmfamily B\kern-.05em \textsc{i\kern-.025em b}\kern-.08em
T\kern-.1667em\lower.7ex\hbox{E}\kern-.125emX}}

\newcommand{\note}[2]{}
\newcommand{\addressed}[3]{}

\def\volumeyear{2024}

\begin{document}

\runninghead{Balos et al.}

\title{SUNDIALS Time Integrators for Exascale Applications with Many Independent ODE Systems}


\author{Cody J. Balos\affilnum{1}, Marc Day\affilnum{3}, Lucas Esclapez\affilnum{3}, Anne M. Felden\affilnum{5}, David J. Gardner\affilnum{1}, Malik Hassanaly\affilnum{3}, Daniel R. Reynolds\affilnum{2}, Jon Rood\affilnum{3}, Jean M. Sexton\affilnum{4}, Nicholas T. Wimer\affilnum{3}, and Carol S. Woodward\affilnum{1}}

\affiliation{\affilnum{1}Center for Applied Scientific Computing, Lawrence Livermore National Laboratory, USA\\
\affilnum{2}Department of Mathematics, Southern Methodist University, USA\\
\affilnum{3}National Renewable Energy Laboratory, USA\\
\affilnum{4}Lawrence Berkeley National Laboratory, USA\\
\affilnum{5}TU Delft, The Netherlands}

\corrauth{Cody Balos,
Lawrence Livermore National Laboratory,
7000 East Ave, L-561, Livermore, CA, USA.}

\email{balos1@llnl.gov}

\begin{abstract}

Many complex systems can be accurately modeled as a set of coupled time-dependent partial differential equations (PDEs).  However, solving such equations can be prohibitively expensive, easily taxing the world’s largest supercomputers. One pragmatic strategy for attacking such problems is to split the PDEs into components that can more easily be solved in isolation. This operator splitting approach is used ubiquitously across scientific domains, and in many cases leads to a set of ordinary differential equations (ODEs) that need to be solved as part of a larger “outer-loop” time-stepping approach. The SUNDIALS library provides a plethora of robust time integration algorithms for solving ODEs, and the U.S. Department of Energy Exascale Computing Project (ECP) has supported its extension to applications on exascale-capable computing hardware.
In this paper, we highlight some SUNDIALS capabilities and its deployment in combustion and cosmology application codes (Pele and Nyx, respectively) where operator splitting gives rise to numerous, small ODE systems that must be solved concurrently.
\end{abstract}

\keywords{Ordinary differential equations, initial value problems, partial differential equations, operator splitting, exascale, GPUs, computational fluid dynamics, chemical kinetics, combustion, cosmology}

\maketitle

\section{Introduction\label{s:intro}}

The accurate simulation of many interesting physical phenomena, such as the evolution of matter in the universe or the combustion of fuel in an engine, can be difficult and costly.
The governing equations for such problems are often in the form of a time-dependent system of partial differential equations (PDEs) that couple multiple physics processes together, where each process may have its own unique and challenging characteristics (e.g., different time scales).
Temporal evolution strategies for these multiphysics problems typically involve one of two distinct approaches.
The first is to apply a time integration method to the full set of spatially discretized equations. The second is to utilize an operator splitting technique \citep{mclachlanSplittingMethods2002} that allows for separate specialized numerical treatments of each physical process.
An operator splitting approach is particularly well suited to computational fluid dynamics (CFD) problems with coupled chemical reactions \citep{earlySplitMeths}.
In these CFD problems, the transport and reactions are split and integrated in time with different numerical methods.
In this case, under standard finite volume or finite difference spatial discretizations, the reactions present as a system of ordinary differential equations (ODEs) at each spatial grid cell of the computational domain,
and their evolution is often the dominant cost of the CFD simulation.
Solving these operator split CFD problems efficiently and at large scales on modern high performance computing systems, in particular systems with both CPUs and GPUs, poses several challenges.
First, a robust and efficient time integrator that can use GPUs is needed for the chemistry evolution.
Choosing an integrator that optimally balances solution accuracy and computing time, is problem dependent \citep{lapointe2020data}, thus the flexibility to try different integration approaches is necessary.
Second, the small size of the numerous independent ODE systems produced by the operator splitting approach does not result in enough work for GPUs on a per unit basis, so a strategy for exposing more concurrency is required.
Thirdly, state information has to be shared efficiently between the transport and reactions.
Within these primary challenges many secondary challenges arise.
For example, when implicit time integration methods are applied, the resulting linear systems have to be solved on the GPU, and either a Jacobian matrix has to be computed or a matrix-free linear solver must be used.




The Pele and Nyx codes are two such applications where the efficient solution of numerous independent ODEs is critical to performance.
The Pele suite of codes \citep{GithubPele} solve reactive flow hydrodynamics problems for combustion applications, while the Nyx code \citep{GithubNyx} solves N-body and gas dynamics problems in cosmology.
The Pele and Nyx codes both employ finite volume methods with adaptive mesh refinement (AMR) to handle spatial discretization of the fluid operators and an operator splitting approach that results in spatially decoupled ODEs in each cell that represent chemical or nucleosynthesis reactions, respectively.
These ODEs must be advanced over a time interval determined by an outer-loop fluid time integration algorithm at every numerical fluid time step during the simulation.
To perform the ODE evolutions Pele and Nyx leverage the SUNDIALS library \citep{hindmarsh2005sundials,balos2021enabling, gardner2022enabling} of solvers which provides access to a variety of a feature-rich time integration methods that efficiently exploit GPUs.
In this paper, we present how the SUNDIALS, Pele, and Nyx teams together developed approaches to tackle a broad range of key challenges that arise when evolving numerous ODE systems.
Specifically, we address the need for flexible, problem-dependent ODE integrators, effective scalable GPU utilization, and efficient communication of state information between the application codes and SUNDIALS. Due to the critical role that SUNDIALS plays in Pele and Nyx, the work discussed here has led to immediate and significant beneficial impact in the applications, in both cases enabling efficient use of exascale class computing hardware.

The close collaboration between the SUNDIALS, Pele, and Nyx teams is part of the
Exascale Computing Project \citep{ECP}, funded by the U.S. Department of Energy,
to accelerate the development of necessary capabilities in scientific application codes, and the numerical software libraries they depend on, especially for simulations on massively parallel computing systems that use hybrid architectures consisting of both CPUs and GPUs.
These hybrid systems comprise the architecture of an increasing portion of supercomputers on the TOP500 list \citep{TOP500},
including the Frontier\footnote{\url{https://www.olcf.ornl.gov/frontier/}} and soon-to-be deployed Aurora\footnote{\url{https://www.anl.gov/aurora}} and El Capitan\footnote{\url{https://asc.llnl.gov/exascale/el-capitan}} exascale supercomputers, for example. 
In such systems, GPUs account for most of the floating point operation performance (90\% or more on many machines) and introduce specific challenges to the efficient concurrent solution of the spatially decoupled ODEs that arise from operator splitting methods, such as those used in Pele and Nyx.

On CPU-based architectures, the solution of the decoupled ODE systems can be ``embarrassingly parallel'' with separate instances of the time integrator applied concurrently to each system in different execution tasks or threads.
However, for GPUs, the mismatch between SIMD execution and the complex logic of advanced time integrators can make these operations difficult to parallelize \citep{stone2013techniques}.
\cite{Niemeyer2014} implemented a time-explicit algorithm for solving the kinetics equations from operator split reacting flow problems, but the performance of their algorithm was a strong function of stiffness in the system.  For non-stiff cases, they observed an order of magnitude acceleration in wall-clock time to integrate large blocks of reacting cells compared to explicit integrations on CPUs. In addition, for moderately-stiff cases their explicit GPU solver was nearly 60$\times$ faster than implicit CPU-based schemes, but with increasing stiffness the observed acceleration dropped dramatically. The decrease in performance was attributed to thread divergence as each cell followed dramatically different paths through time. These issues, as well as the practical complexity of porting a more general integrator software package to GPUs, suggest an entirely different solution strategy, where the uncoupled systems are batched together.
Batching, in-turn, leads to new challenges for the numerical methods as well as the parallel execution that must be addressed for robust, accurate, and efficient time integration.

In this article, we present a batched strategy for integrating the ODEs that result from an operator split approach to solving coupled reacting flow PDEs. Our approach is tailored to modern GPU-based architectures and presented in the context of applications in cosmology and combustion. The goal of this work is to highlight challenges and key innovations required
to effectively exploit GPU architectures for this class of problems on exascale systems.
Additionally, this work provides a reference for codes that face similar challenges, and, to that end, this paper is structured with elements of a case study.
Section \ref{s:sundials_strategy} describes the overall strategy for solving many independent ODEs from the time integrator point-of-view.
Section \ref{s:analyzing} introduces the Pele and Nyx codes in more detail and describes the distinct challenges in integrating their ODE systems on modern high performance computing hardware.
Section \ref{s:improving} then discusses innovations that improve integrator performance in Pele and Nyx.
The cumulative effects of these improvements are demonstrated in Section \ref{s:results} with representative large-scale runs of the Pele and Nyx codes on the Frontier exascale supercomputer.
Finally, in Section \ref{s:conclusion} concluding thoughts are discussed along with future work.

\section{Overview of SUNDIALS Strategy for Independent ODE Systems on GPUs}\label{s:sundials_strategy}


The SUNDIALS library consists of six packages: CVODE, CVODES, IDA, IDAS, ARKODE, and KINSOL.
CVODE and ARKODE solve ODEs, IDA solves differential algebraic equations, and KINSOL solves systems of nonlinear algebraic equations.
CVODES and IDAS are variants of CVODE and IDA that include sensitivity analysis capabilities.
As this work is concerned with the forward solution of ODEs, we focus here on the methods provided by CVODE and ARKODE.

CVODE solves stiff and nonstiff problems in the form
\begin{equation}
    \frac{\mathrm dy}{\mathrm dt} = f(t, y),\quad y(t_0) = y_0
\end{equation}
where $y \in \mathbb{R}^N$ and $f: \mathbb{R} \times \mathbb{R}^N \to \mathbb{R}^N$
with adaptive order and step size linear multistep methods.
ARKODE is designed for ODEs in the linearly implicit form
\begin{equation}
    M(t)\frac{\mathrm dy}{\mathrm dt} = f(t, y),\quad y(t_0) = y_0.
\end{equation}
and provides adaptive step size diagonally implicit, explicit, additive implicit-explicit,
and multirate Runge--Kutta methods.

Both CVODE and ARKODE base their time step adaptivity on estimates of the local truncation error (LTE) and the user-defined
integrator tolerances \citep{hindmarsh2005sundials,reynolds2023arkode}.
The relative tolerance, $rtol$, controls the LTE relative to the size of the solution.
The absolute tolerance is defined per component of $y$ so that $atol_i$ controls the LTE for the $i$-th solution component when that component is small; $atol_i$ effectively sets the floor below which the solution component will not be resolved.
The integrators will accept a time-step if $\| LTE \|_{wrms} < 1$ where the weighted root-mean square (WRMS) norm is given by
\begin{equation}\label{eq:wrms}
\begin{gathered}
    \| v \|_{wrms} = \sqrt{ \frac{1}{N} \sum_{i=1}^N (w_i v_i)^2 },\\
    w_i = \frac{1}{rtol |y^{n-1}_i| + atol_i},
\end{gathered}
\end{equation}
where $N$ is the number of elements in the vector and $y^{n-1}$ is the previous time step solution. For the remainder of the paper, we will use $\| \cdot \|$ to denote the WRMS norm. If $\| LTE \| \ge 1$, then the step is rejected, the step size is cut (and with CVODE the method order may be reduced), and the step is reattempted. If the step is accepted, the next time step size (and the method order with CVODE) is selected to try to satisfy the error test condition and minimize computational cost.
In most applications, especially computational chemistry, proper integrator tolerances can be the difference between excellent and unacceptable performance.
We demonstrate this in Section \ref{ss:improving_pele}.

With implicit methods, CVODE and ARKODE must solve a nonlinear system, $F(y^n) = 0$, at least once per time step for the future time solution, $y^n$.  For stiff systems, Newton's method is typically the preferred nonlinear solver as it gives much faster convergence than stationary iterations on such problems. The stopping criteria for Newton's method is based on the user-supplied relative and absolute tolerances, to avoid extraneous computational effort while preventing nonlinear solver errors from polluting the temporal error estimates.
Specifically, at the m-$th$ Newton iteration, the stopping criteria is
\begin{equation}\label{eq:stop_newton}
\begin{gathered}
    R \| \delta_m \| < c_{\epsilon} \epsilon,\\
    R = \max\{c_r R, \| \delta_m \| / \| \delta_{m-1} \|\},
\end{gathered}
\end{equation}
where $\delta_m = y^{n(m)} - y^{n(m-1)}$, $\epsilon$ is a method-dependent constant, and the factors $c_\epsilon$ and $c_r$ are constants that default to $0.1$ and $0.3$, respectively.
Each iteration of Newton's method requires the solution of a linear system of equations.
When an iterative linear solver is employed, the integrator tolerances are used to construct scaling vectors to balance the error between solution components within the linear solve and improve efficiency.
Thus, instead of solving $Ax = b$, the iterative linear solvers in SUNDIALS consider the system $\tilde{A} \tilde{x} = \tilde{b}$ with
\begin{equation}\label{eq:linear}
    \tilde{A} = S_1 P_1^{-1} A P_2^{-1} S_2^{-1},\quad \tilde{b} = S_1 P_1^{-1} b,\quad \tilde{x} = S_2 P_2 x
\end{equation}
where $S_1$ and $S_2$ are diagonal scaling matrices with entries $w_i$, and $P_1$ and $P_2$ are (optional) application-specific left and right preconditioners, respectively.
The iterative linear solver stopping criteria for this modified system is
\begin{equation}\label{eq:stop_linear}
\begin{split}
    \|\tilde{b} - \tilde{A}\tilde{x}\|_2 < c_{l} (c_\epsilon \epsilon)
\end{split}
\end{equation}
where $c_l$ is a constant (defaults to $0.05$) that scales the nonlinear solver stopping tolerance from Equation~\ref{eq:stop_newton} to ensure sufficient accuracy without over solving the linear system. Note the SUNDIALS GMRES implementation computes the linear residual norm in Equation~\ref{eq:stop_linear} from the rotations used to solve the least-squares problem \citep{saad2003iterative} and thus avoids additional reduction operations.
As we show in Section \ref{ss:improving_pele}, this linear solver scaling strategy is critical in many application settings, including Pele, where the various components of the solution can have dramatically disparate ranges of scale.

As is the case with all of the SUNDIALS packages, CVODE and ARKODE are based on a set of shared SUNDIALS base classes that provide utilities (e.g., memory management), data structures such as vectors and matrices, and both nonlinear and linear algebraic solvers \citep{gardner2022enabling}.
This design encapsulates both the operations on data and the parallelization employed by the time integration and solver algorithms. Modification of state data, $y$, occurs within class implementations or in application provided callback functions to evaluate, for example, the ODE right-hand side, $f$,  or Jacobian, $\partial f / \partial y$.
In the context of a CPU-GPU system, SUNDIALS is designed so that the complex integrator logic is executed on the CPU, but GPU kernels are launched (through the class implementations or callback functions) to operate on state data that lives in GPU memory \citep{balos2021enabling}.
There are a couple of reasons this execution model was chosen.
Foremost, thread divergence caused by the complex timestep adaptivity and solver heuristics within the implicit SUNDIALS integrators limits the efficiency of the
alternative execution model discussed in Section~\ref{s:intro}, where the entire integrator and solver stack executes on the GPU \citep{stone2013techniques}.
The secondary reason was to ensure sustainability and maintainability of SUNDIALS. With the chosen execution model (integrator logic on CPU, data operations on GPU) the tens of thousands of lines of code in the integrators do not have to be re-implemented as GPU kernels, and long-standing SUNDIALS data structures can be leveraged.

The design of SUNDIALS has several implications when used with hybrid CPU-GPU systems.  For the applications considered here, the time splitting strategies require that one ODE system be advanced for each cell in the domain, over a global time interval, $dt_{\textit{CFD}}$ (determined by the transport physics).  This task is divided into work across batches of cells where each batch is advanced independently by SUNDIALS. The number of cells in each batch is determined by the application, but should be large enough to guarantee that computational work on each GPU dwarfs the kernel launch overheads yet is small enough that each batch advance fits entirely into the available GPU memory. Within each batch, all cells are integrated in lock-step with subcycled time steps, $dt_{\textit{chem}}$, that are selected adaptively by the SUNDIALS algorithms as the integration over $dt_{\textit{CFD}}$ proceeds. At any instant, $dt_{\textit{chem}}$ is typically constrained by the stiffest cell in the batch. In cases where the stiffness varies widely over the domain, application-specific batching and load balancing strategies should be used to evenly distribute the total integration work.
Note also that as a consequence of batching ODE systems, the SUNDIALS approach to parallelism on CPU-GPU systems suggests that special linear solvers that exploit the batched structure can be critical for computational efficiency.

\section{Analyzing the Applications}\label{s:analyzing}


We now discuss the Pele and Nyx applications and their uses of SUNDIALS for solving many independent ODE problems.
Many decisions regarding algorithms and parameters arise when solving the ODEs in these applications.
For the remainder of this paper we will refer to a set of algorithms and parameter choices as a solution approach.
In particular, we are concerned with the time integration method, nonlinear solver, linear solver, and tolerance selection (which is somewhat unique to Pele).
Table \ref{tab:sol_app} lists the primary solution approaches examined.

\begin{table*}[!htb]
\centering
\caption{\label{tab:sol_app} Implicit approaches to solving the ODE systems in Pele and Nyx require a time integration scheme, nonlinear and linear solvers, a way to compute the Jacobian (or its action on a vector), and a strategy for selecting the integrator tolerances (discussed in Sec.~\ref{sec:tolerances}). Here we present the different solution approaches discussed in this paper.
}
\begin{tabular}{l|l|l|l|l|l}
Solution Approach&Integrator&Nonlinear Solver&Linear Solver&Jacobian&Tolerance Strategy\\
\hline
1A&CVODE BDF&Inexact Newton &GMRES&Numerical Jv-product&Fixed\\
1B&CVODE BDF&Inexact Newton &GMRES&Numerical Jv-product&Typical Values\\
2A&CVODE BDF&Modified Newton&Dense Direct&Analytical&Fixed\\
2B&CVODE BDF&Modified Newton&Dense Direct&Analytical&Typical Values\\
3A&CVODE BDF&Modified Newton&Dense Direct&Numerical&Fixed\\
3B&CVODE BDF&Modified Newton&Dense Direct&Numerical&Typical Values\\
\end{tabular}
\end{table*}

\subsection{Pele}
\label{ss:analyzing_pele}
The Pele Suite is a collection of codes and algorithms for simulation and analysis of turbulent reacting flows
and has been co-designed throughout to take advantage of modern exascale computing hardware. The core capabilities in the Pele Suite are two block structured adaptive mesh flow solvers; PeleC~\citep{PeleC_IJHPCA} for compressible flows and PeleLMeX~\citep{Esclapez2023} for low Mach systems. Both PeleC and PeleLMeX leverage the AMReX library \citep{Zhang2019} to provide distributed data containers, associated algorithmic infrastructure, and performance portability for efficient execution on both CPU and GPU hardware. 

PeleC solves the fully compressible, multi-species reacting Navier-Stokes equations for the conservation of species mass, momentum, and energy, and includes finite rate chemistry. A finite-volume discretization is used for the spatial operators and the time advance is based on an operator splitting strategy. Advection and diffusion terms are treated explicitly in time, while chemical reaction terms are integrated over each flow solver time step, $dt_{\textit{CFD}}$, using ODE solvers from SUNDIALS, or a built-in explicit integrator when the stiffness is limited.
An outer loop iteration is applied to couple the flow and chemistry integration with a scheme that is second-order in space and time.

In contrast, PeleLMeX solves the reacting flow equations in the low Mach number limit where acoustic waves are assumed to traverse the domain infinitely fast relative to transport processes. Because sound waves are eliminated from the solution, the fluid solver time steps, $dt_{\textit{CFD}}$, taken in PeleLMeX can be considerably larger than those in PeleC. An iterated variable-density projection scheme is used in PeleLMeX to advance a second-order finite-volume discretization subject to the elliptic constraint of spatially constant thermodynamic pressure. Diffusion is treated implicitly so that $dt_{\textit{CFD}}$ is limited only by the explicit treatment of advection. Similar to PeleC, ODE integrators from SUNDIALS are used in PeleLMeX to integrate the chemistry with a temporal discretization that is operationally split from the flow advance. Similar to the approach described in
\citep{Zingale2022}, an outer iteration ensures that each component of the operator split treatment is sufficiently updated to guarantee second-order convergence of the overall advance. Importantly, integrating the chemistry in PeleLMeX can be considerably more taxing computationally than in PeleC due to the longer integration times resulting from the larger fluid solver time steps.

In both Pele codes,
the model equations can be considered as a general additively partitioned system
\begin{equation}
\frac{\partial U}{\partial t} = F + R,
\label{chemODE}
\end{equation}
where $U$ is the state vector, $R$ represents the chemical source terms, and $F$ represents all forcing on the state other than chemistry (i.e., advection, diffusion, etc.). When discretizing Equation~\ref{chemODE}, $F$ depends only on the current state in PeleC while PeleLMeX utilizes a semi-implicit approach where $F$ is computed with a lagged $R$ and updated iteratively.  That is, a provisional $F$ is computed with an estimate of $R$, the chemical source term averaged over $dt_{\textit{CFD}}$. Over this interval, $R$ is recomputed given the new $F$ as an averaged source term, and the process is iterated until both terms cease to change up to a user-specified tolerance.
While computing $R$, $F$ is time-centered and completely known (and treated as temporally constant) at each cell in the domain during the ODE solution over $dt_{\textit{CFD}}$ using SUNDIALS. Note specifically that any stencil operations used to compute $F$ have already been reduced to incorporate any dependence on neighbor cells. For both PeleC and PeleLMeX, Equation~\ref{chemODE} is thus a spatially decoupled ODE.

In many combustion problems, $R$\, in Equation~\ref{chemODE} exhibits a broad range of timescales and can often be numerically stiff in some parts of the domain (e.g., at a flame surface). Due to stiffness, non-linearities and accuracy\slash realizability constraints, such as the need for positive species compositions, numerical integration of this ODE often consumes more than 70\% of the total resources required for the simulation.

In addition to routines that evaluate $F$ and $R$ in Equation~\ref{chemODE}, some of the ODE integration methods available in SUNDIALS require additional information about the system, including the Jacobian matrix $\partial R/\partial U$ and an estimate for the magnitude of each component $i$ of $U$ (used to determine $atol_i$ in Section \ref{s:sundials_strategy}). Each of these requirements also suggests opportunities for optimizations to accelerate the ODE solution process.

Because the ODEs are uncoupled across the domain, we have considerable freedom in staging their integration over $dt_{\textit{CFD}}$. In multi-threaded CPU-based implementations, the total work to integrate all the cells in the domain can be distributed arbitrarily across the available threads with no race or synchronization concerns. The work for each cell can be approximated ahead of time by tracking the work required to integrate the previous CFD time step. Thus, load balance in the CPU implementation can be achieved by evenly distributing the predicted work across processors, and within each processor by distributing across available threads.  As we indicate below, the situation is more complex in the GPU implementation.


For the CPU implementation, other simple optimizations can have considerable impact on simulation run times.  With relatively limited parallelism, each thread can be assigned to integrate a set of cells sequentially. Auxiliary space may be required to compute and temporarily store the system Jacobian. Formally, this Jacobian should be recomputed for every Newton iteration of every subcycled time step taken by the ODE integrator, but since only one cell is integrated at a time, this space can be allocated for a single cell and reused for all subsequent cells integrated by each thread. Also, because computing the Jacobian can be expensive, a standard optimization is to apply a modified Newton method which attempts to reuse the Jacobian for multiple Newton iterations, and even over multiple subcycled time steps, unless that reuse results in local divergence of the nonlinear iterations.
This approach is the default implementation of Newton's method when using a matrix-based linear solver in CVODE and ARKODE.
In many situations, however, a sequence of cells in the domain are similar enough and the ODE is ``easy'' enough that a further optimization can be done by using the same Jacobian across a large number of cells. This simple strategy can reduce the ODE integration effort by factors of 4-10 or more because many combustion problems have large regions of nearly uniform composition and temperature (e.g., cold fuel, hot products, etc.). Similarly, parts of the evaluation of $R$ and the Jacobian often involve Arrhenius expressions requiring expensive transcendental function evaluations that also can be cached and reused across adjacent cells. In practice, this reuse can reduce compute efforts by another factor of 2-3.

Finally, in combustion problems, different components of the state, $U$, can take on radically disparate scales --  temperature is typically 300 to 2500 K, whereas species mass fractions could span values of 1 to 10$^{-13}$.  Without proper scaling of the various norms and tolerances used within the ODE integrator, one would need to converge all components to the same small tolerance.
Satisfying this requirement could be excessively wasteful computationally or even impossible due to numerical roundoff. This scaling is affected through setting of appropriate scale values feeding into the absolute tolerance values in
Equation~\ref{eq:wrms}.
As noted above, the scale values are used to scale all integrator error norms, to compute algebraic solver tolerances, to compute perturbations used for directional derivatives by the Newton-Krylov methods, and to compute numerical Jacobian entries, if applicable.
Setting proper values for these scale factors can have a massive impact on the robustness and efficiency of the ODE integration, improving performance by factors of 100$\times$ or greater.
Note, however, that ``proper'' scale values may vary by orders of magnitude over the domain and over time, even in a single simulation. In Section~\ref{s:improving} we provide a simple strategy that can be used to approximate and periodically update estimates of the scale of each state component based on the evolving global solution.

As discussed in Section~\ref{s:sundials_strategy}, the batched system strategy used with SUNDIALS and GPU systems is designed to maximally exploit the inherent parallelism of the device.
Rather than thread-based groups of cells being integrated sequentially, parallelism across cells is mapped to the individual GPU cores.
That is, each cell is integrated by a single GPU core. This strategy has a number of implications on the optimizations mentioned above as well as those within the SUNDIALS integration algorithms.  Firstly, any strategy depending on the caching and reuse of expensive computation no longer applies since all cells within each batch are integrated in lockstep.  The simple load balancing strategy based on cell-wise work estimates no longer applies directly either; neither do published strategies for spatial variable scaling. Storage for system Jacobians, if required, is needed across all cells at once, and this dramatically increases the memory required to work on a given batch of cells. Finally, since the integration is lockstep across a batch of cells, the effort required for numerical subcycling, and even local nonlinear solves will be dictated by the most difficult cells in the batch. Of the simple optimization and load-balancing methods discussed above, only variable scaling remains as a strategy that we have employed in this work.
We discuss examples of this scaling in subsequent sections.


\subsection{Nyx}
\label{ss:analyzing_nyx}


Nyx is a highly parallel, adaptive mesh, finite-volume N-body compressible hydrodynamics
solver for cosmological simulations. It has been used to simulate different cosmological scenarios with a recent focus on the intergalactic medium \citep{Onorbe2017} and Lyman-$\alpha$ forests \citep{Almgren2013,Sexton2021}. Similar to the Pele Suite of codes, Nyx leverages the AMReX library for CPU and GPU performance-portability.  However, Nyx makes more extensive use of AMReX's particle\slash mesh capabilities in order to model the interaction of gas dynamics with cosmological dark matter.

Nyx solves the compressible hydrodynamic equations on an adaptive grid hierarchy coupled with an N-body treatment of dark matter. The gas dynamics in Nyx use a finite volume methodology on a set of 3-D Eulerian grids based on the directionally unsplit corner transport upwind method of Colella with piecewise parabolic reconstruction \citep{ppm}. Dark matter is represented as discrete particles moving under the influence of gravity and evolved via a particle-mesh method, using a Cloud-in-Cell deposition/interpolation scheme.
In addition, Nyx includes physics needed to accurately model the intergalactic medium in the optically thin limit and assuming ionization equilibrium.
These simplifying assumptions reduce the model for heating-cooling to a single scalar equation for each physical cell location, which is solved using SUNDIALS and coupled to the hydrodynamics using a modified spectral deferred corrections (SDC) scheme \citep{zingale2019}.



Similar to the Pele codes, the coupled hydrodynamic and heating-cooling processes in Nyx are treated as an additively partitioned system
\begin{equation}
\frac{\partial e}{\partial t} = F + R
\label{nyxODE}
\end{equation}
advanced over a time interval determined by the discretized physics represented in $F$. Here $e$ is the internal energy, $R$ is the heating-cooling source terms, and $F$ represents all other forcing on the state, e.g., hydrodynamics. In Nyx, $F$ is computed with a semi-implicit time discretization over the interval, $dt_{\textit{CFD}}$.  A lagged approximation of $R$ (averaged over $dt_{\textit{CFD}}$) is added as an external forcing for the $F$ update. The heating-cooling solve then proceeds with this approximation to $F$ as an external forcing.  Note specifically that for the ODE solve, $F$ is thus completely known, time-centered, and treated as temporally constant at each cell. Additionally, since any stencil operations have already been incorporated into the computation of $F$, Equation~\ref{nyxODE} is a spatially decoupled ODE.

In many cosmological problems, $R$\, in Equation~\ref{chemODE} exhibits a broad range of timescales and can often be numerically stiff in localized regions of space.  Moreover, the stiffness can vary dramatically over time, such as in cooling condensing systems, particularly when the cooling models have strong density and redshift dependencies.
Due to stiffness, nonlinearities, accuracy and realizability constraints, such as the need for non-negative temperatures and densities, numerical integration of the spatially decoupled ODEs representing heat-cooling processes can often consume as much as 30\% of the total resources required for a Nyx simulation.

\subsection{Common Themes and Differences}
\label{ss:analyzing_common}

As discussed, the Pele and Nyx codes share a common overall code structure and a similar operator split strategy for incorporating reaction terms into their time integration algorithms.  Since they both are built upon AMReX, they also share common underlying data structures used to define and manage the state.  Both codes store the state in logically rectangular containers corresponding to a uniform Cartesian grid, with an underlying 1D array structure in memory that maps onto component and index space. By default in AMReX, the state component varies slowest in this mapping in order to optimize local stencil operations. In order to maximize cache reuse during the reaction integration the first step in both codes is to create and fill an auxiliary (temporary) memory block with reordered state data such that all components at a cell are contiguous in memory. Due to the potentially large number of species in the Pele state vector, temporary memory associated with this reordering can substantially impact the maximum number of cells that can be stored on each core. Future research should explore the tradeoff in memory available for this operation and the gains due to improved cache reuse that results.  In both Pele and Nyx, the stiffness and/or computational difficulty of integrating each cell can vary dramatically over the domain. However, since the ODE solves are batched and integrated in lockstep over large numbers of cells in parallel, the overall effort to advance a block of cells is determined by the most difficult ones in the set. Notwithstanding performance improvements discussed in the next section, future research could also explore various strategies to batch cells expected to have similar dynamics and stage the integration of these batches for maximum overall performance.

\section{Improving Performance}\label{s:improving}

In this section we discuss some of the algorithm and software changes implemented to improve the performance of the ODE integrations in Pele and Nyx.

\subsection{Improving performance in Pele}
\label{ss:improving_pele}



\subsubsection{Load balancing}
As mentioned previously, load balancing is not handled internally by the SUNDIALS integrators but rather is the responsibility of the application developer. In the following, we present an example load balancing strategy used in PeleLMeX, for illustration purposes. First, it is useful to summarize a typical hierarchy of spatial mesh refinement levels.

The AMR grid structure is built as a set of levels. Each level consists of a union of non-overlapping boxes that each define a rectangular subregion of index space. State data is stored at each index in these boxes. The data is distributed across the parallel machine at the granularity of a box, and is mapped to a corresponding physical location assuming constant cell spacing. The coarsest AMR level covers the entire computational domain. The next finest AMR level also consists of a union of boxes, but with finer grid spacing, and it typically covers only a subset of the domain, focused around an interesting (possibly dynamic) feature in the solution. Data in the boxes at this level is also distributed by box across the machine.  Increasingly finer levels cover increasingly smaller regions of the domain until all interesting features are covered by appropriately resolved grids. Periodically, the level structure is reconstructed in order to maintain the desired resolution of the various evolving solution features. Note that the finest AMR grid patches are typically focused around regions where the computational burden of the ODE integrator is the highest. The size (typically 8 to 32 cells on a side), number, and distribution of these grid patches are determined by the AMReX library based on user-defined parameters \citep{Zhang2019}. On GPU platforms, larger grid patches are favored to exploit the GPU's high degree of parallelism.

In Pele, for all stencil-based operations (advection, diffusion, etc.) these grid patches are usually distributed across processors using a space curve filling approach in order to minimize inter-processor communications.
For the ODE chemistry integration, a second auxiliary set of grid patches (called the ``ODE-grid''), covering the same domain real estate, is introduced and distributed across processors to account for the following:
\begin{itemize}
    \item AMR level hierarchies are advanced in time with either an assumption of constant CFL or constant $dt_{\textit{CFD}}$.  The former requires ``subcycling'', where finer grids are advanced with smaller time steps and physics-specific adjustments are orchestrated to ensure satisfaction of global constraints (e.g., conservation). Alternatively, all levels are advanced together and cells in the coarser levels covered by fine cells are irrelevant to the solution. These ``covered coarse'' cells can be removed from the CFD and chemistry advance.
    \item The ODE integration is fully local such that communication is not a concern allowing for other approaches for distributing work.
    \item The estimated ODE integration work of each grid patch is not solely proportional to its number of cells, but also includes an estimate of the integrator effort (e.g., the number of calls to the chemistry right-hand-side function in the previous step);
    \item The patch sizes in the ODE-grid are selected to enable control over the granularity of the load balancing and to help mitigate the memory requirements of the ODE solver while also maintaining enough workload for GPUs.
\end{itemize}
The first three points are relevant to both CPU and GPU simulations, but the last constraint is specific to GPUs.
On CPUs, a workload estimate is available for each grid cell where an independent ODE is solved, allowing fine-grained control of the load distribution when building the ODE-grid with grid patches as small as possible. On GPUs, an ODE is solved for all the cells in a given patch, reducing the granularity of the work estimate and, thus, the potential gain from improved load balancing. Reducing the ODE-grid patch size to gain more granularity on GPUs eventually degrades the performance, since this decreases the parallelism exposed to the GPU and increases the relative cost from kernel launch overhead. The size of the ODE-grid patches is a user-defined parameter that depends on the size of the chemical mechanism, the linear solver choice in the ODE integrator, and other case-specific parameters. In practice, the ODE-grid is updated at every AMR regrid operation and is revisited at regular intervals.  Here, a single processor will build a tentative ODE-grid and trigger the deployment of a new ODE-grid if the tentative ODE-grid improves the load balancing enough (typically, by 5\% or higher to balance the cost of data movement incurred by deploying the new ODE-grid). We note that this strategy alone provided a 10\% to 15\% reduction of the overall PeleLMeX compute time while running an Exascale-relevant auto-ignition case on Frontier.


\subsubsection{Data Ordering}
\begin{figure}[!htb]
\centering
\includegraphics[width=0.45\textwidth]{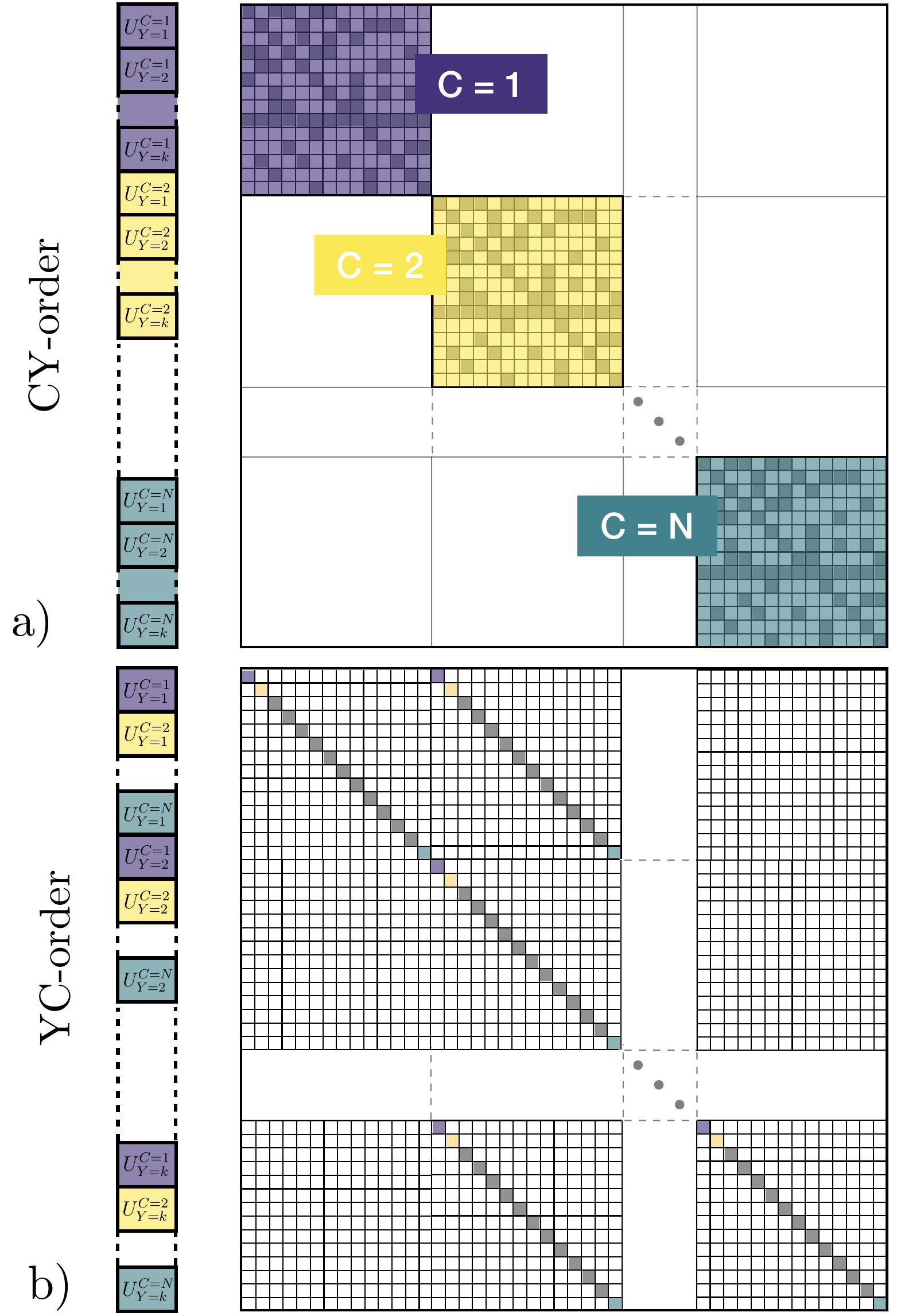}
\caption{a) CY-ordering for the state vector, where all the state components of a given cell are adjacent in the state vector, along with a schematic of the Jacobian pattern. b) YC-ordering, where adjacent entries contain a given state component for all the cells followed by another component
for all the cells, along with the Jacobian pattern.
}
\label{fig:pele_ordering_PP}
\end{figure}

SUNDIALS integrators effectively interweave calls to reduction operations on the vector data structure and calls to the Pele system right-hand-side (RHS) or Jacobian functions. On a GPU, multiple individual computational cells are bundled within a single ARKODE or CVODE instance, so the data layout should be adapted to enable coalescing memory access in the application kernels when possible.
In particular, Pele kernels use a cell-wise granularity with each GPU thread handling the RHS kernel for the entire chemical system in each cell. To allow coalescing memory access, the state vector must then be ordered such that adjacent entries contain a given state component for all the cells followed by another component for all the cells and so on (\emph{YC}-order). However, as illustrated in Figure \ref{fig:pele_ordering_PP}(b), this ordering is incompatible with the use of so-called ``batched linear solvers'' within the modified Newton step, that expect repeated Jacobian blocks having the same nonzero structure.  In contrast, having all the state components for a given cell adjacent in the state vector (\emph{CY}-order) results in the desired block-diagonal structure, as illustrated in Figure \ref{fig:pele_ordering_PP}(a). Both of these orderings are available in Pele; \emph{YC}-order is the default for matrix-free solution approaches (1A and 1B from Table \ref{tab:sol_app}), and \emph{CY}-order is preferred when using solution methods that rely on assembly of Jacobian matrices (2A, 2B, 3A, and 3B from Table \ref{tab:sol_app}).

\subsubsection{Selecting Integrator Tolerances}
\label{sec:tolerances}

\begin{figure*}[!htb]
\centering
\includegraphics[width=0.49\textwidth]{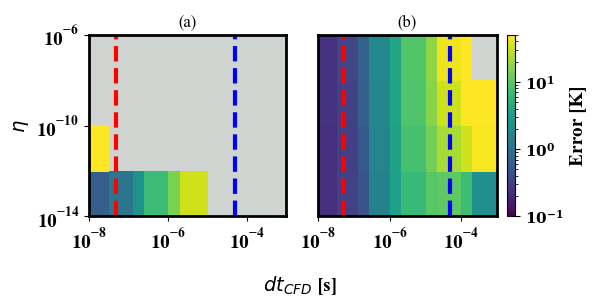}
\includegraphics[width=0.49\textwidth]{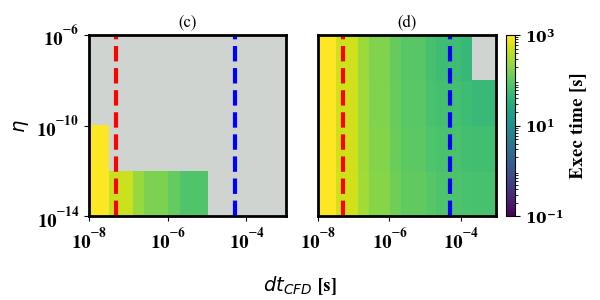}
\caption{Dynamically setting the CVODE absolute tolerances by sampling typical values in Pele results in better accuracy and is faster than using a fixed value for the tolerance in nearly every scenario tested.
The contour plots (a) and (b) show the average mean-square error of the zero-dimensional reactor temperature profile against $\eta$ (see Equation~\ref{eq:eta}) and the fluid timestep, $dt_{\textit{CFD}}$, for solution approaches 1A and 1B from Table \ref{tab:sol_app}. The contour plots (c) and (d) show the execution time. Gray indicates that the run timed out, or that the error was so large that there was no ignition of the reactants. The typical $dt_{\textit{CFD}}$ for PeleC (red dashed line) and PeleLMeX (blue dashed line) are shown for reference.
}
\label{fig:err_sol}
\end{figure*}

\begin{figure}[!htb]
\centering
\includegraphics[width=0.45\textwidth]{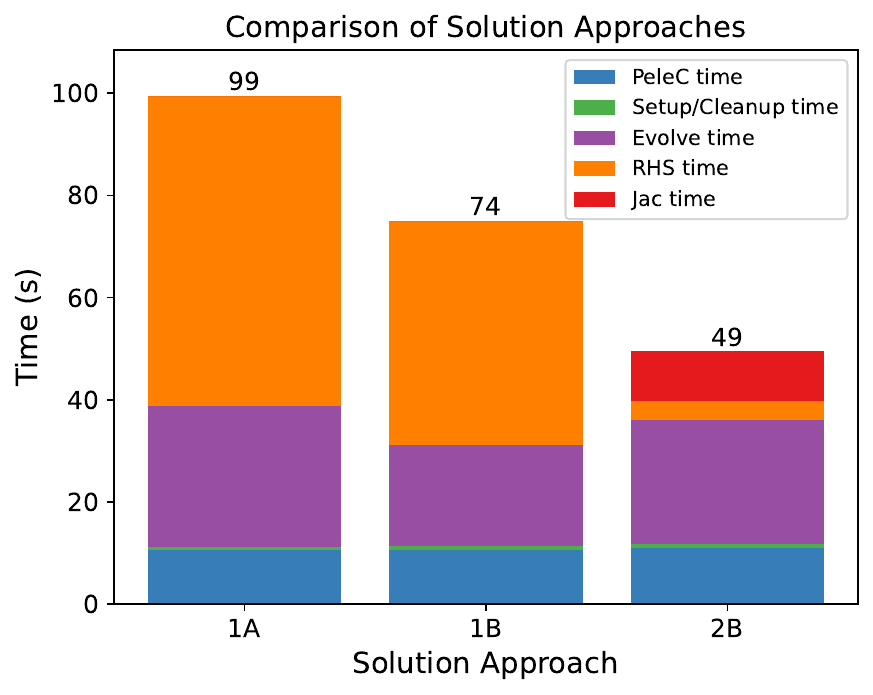}
\caption{PeleC timings with CVODE on a pre-mixed flame test problem on AMD MI100 GPUs using three of the solution approaches from Table~\ref{tab:sol_app}. Approach 1A uses GMRES with fixed absolute tolerances, approach 1B also uses GMRES but additionally uses typical value absolute tolerances, and 2B uses a direct solver (MAGMA) and typical value absolute tolerance. All three approaches use shared memory reductions. Using typical value absolute tolerances improves the performance by 25\% while the direct solver gives the fastest results.}
\label{fig:pelec_cvode_compare}
\end{figure}


As alluded to in Section \ref{ss:analyzing_pele}, selection of appropriate ODE integrator tolerances is critical in light of the disparity of scales between the temperature and individual species mass fractions. Choosing a single scalar value for all the state entries appearing in Equation \ref{eq:wrms} leads to improper scaling of norms both in the iterative linear solvers and in the integrator error test, possibly leading the ODE integrator to use smaller time steps than necessary to achieve the desired accuracy.
To avoid this small step size restriction, Pele defines individual absolute tolerance values for each state component. These tolerances rely on the so-called ``typical values'' of the individual chemical species or temperature variable, denoted $\tilde{y}_i$, where
\begin{equation*}
    \tilde{y}_i = \frac{1}{2}\big( \min(y_i) + \max(y_i) \big),
\end{equation*}
and the $\min$ and $\max$ operations are taken over the entire computational domain.
The typical values are reset at user-defined intervals.
The absolute tolerance is then set using a scaling factor, $\eta$, as
\begin{equation}
    atol_i = \eta\, \tilde{y}_i.\label{eq:eta}
\end{equation}
The value of $\eta$ is problem dependent (representative values are given in the examples shown later), but
the default choice adopted in Pele is $\eta = 10^{-10}$. When typical values are not explicitly managed and provided, then $\tilde{y}_i$ defaults to 1 in Equation~\ref{eq:eta}.  Note that between updates to the typical values, the $atol_i$ remain constant as the ODEs evolve within SUNDIALS.


To demonstrate the importance of the typical value tolerances on CVODE performance in Pele, we simulated a zero-dimensional reactor with a 53 species $n$-dodecane mechanism \citep{yao2017compact} for $t \in [0, 0.1]$s, initialized with a mixture at an equivalence ratio of 7.93, initial temperature of 600K and initial pressure of 5 atm. The simulations assumed constant internal energy.
This reactor simulation is similar to what would arise during a full PeleC or PeleLMeX simulation.
To emulate the different fluid time steps, $dt_{\textit{CFD}}$, taken by PeleC and PeleLMeX the time interval is subdivided into uniform intervals ranging from 10$^{-8}$s to 10$^{-3}$s and the subintervals are solved sequentially by CVODE. Specifically, we emulate $dt_{\textit{CFD}} = \{$1, 3.125, 12.5, 25, 100, 200, 1000, 20,000, 100,000$\} \times 10^{-6}$s.
The test is performed for a range of absolute tolerance scale factor levels ($\eta$ in Equation~\ref{eq:eta}) taken from the following set $\{10^{-6}, 10^{-8}, 10^{-10}, 10^{-12}, 10^{-14}\}$.  The relative tolerance is held fixed at $10^{-7}$ for all tests.
We employed two different solution approaches and evaluated the impact of the tolerances in terms of accuracy and execution time:

\begin{enumerate}
    \item CVODE BDF with an inexact Newton-Krylov nonlinear solver with a fixed scalar value $atol_i = \eta$ (i.e., $\tilde{y}_i=1$, solution approach 1A in Table \ref{tab:sol_app}).
    \item CVODE BDF with an inexact Newton-Krylov nonlinear solver with dynamic specification of $atol_i$ for each state component via Equation~\ref{eq:eta} using the typical values (solution approach 1B in Table \ref{tab:sol_app}). In this case, the typical values are computed based on the state at the end of the reference integration.
\end{enumerate}

The accuracy of the simulations is evaluated by comparing the error in temperature profiles over the integration time (0.1s) to a reference simulation.
We define the error as the averaged mean squared error of temperature over all the time steps simulated. It can be understood as the average temperature error incurred at each fluid cell in a multidimensional simulation.
The reference simulation used the modified Newton solver in CVODE with the Jacobian computed using local finite differences (solution approach 3A from Table \ref{tab:sol_app}). To ensure vanishing numerical errors in this reference solution, it is integrated using atypically small uniform intervals of $dt_{\textit{CFD}}=3.3\times 10^{-9}$s and the absolute and relative tolerances are set to $10^{-15}$ for all solution components.
Test configurations which result in unacceptable execution times (over 7200 seconds, arbitrarily) are terminated. Integration ``failure'' was also declared if the errors were so large that there was no ignition of the reactants over the entire interval.

Over the $dt_{\textit{CFD}}$ and tolerance ranges chosen, the error can vary by several orders of magnitude (Figure~\ref{fig:err_sol} (a) and (b)).
Likewise, the execution time for each case depends on the tolerance and the $dt_{\textit{CFD}}$ used (Figure~\ref{fig:err_sol} (c) and (d)). Unsurprisingly, the most computationally intensive settings are also the most accurate.
Furthermore, our typical value tolerance strategy results in faster execution times while maintaining accuracy for all of the $dt_{\textit{CFD}}$ and $\eta$ values except when $\eta$ is very small.
The effect of $dt_{\textit{CFD}}$ and tolerance level on the execution time of the chemistry integration has also been observed in other reacting flow solvers \citep{lapointe2020data}.

In the context of a 3D PeleC simulation of a simple harmonically perturbed pre-mixed flame, the typical value tolerances reduce the run time by approximately 25\% when using the Newton-Krylov method in SUNDIALS and GMRES as the linear solver (Figure~\ref{fig:pelec_cvode_compare}).
The reduction in execution time is due primarily to an approximately 30\% reduction in the average number of GMRES iterations per ODE solve with CVODE.
The improved linear solver performance also results in fewer algebraic solver convergence failures and thus fewer failed steps and RHS evaluations.
The impact of the tolerance selection on the number of iterative linear solver iterations is due to WRMS norm used in the nonlinear solver stopping criteria (Equation~\ref{eq:stop_newton}) and the scaling matrices applied to the linear system (Equation~\ref{eq:stop_linear}).

\subsubsection{Solver Robustness}
\label{sec:solv_rob}

\begin{figure}[!htb]
    \centering
    \includegraphics[width=0.47\textwidth]{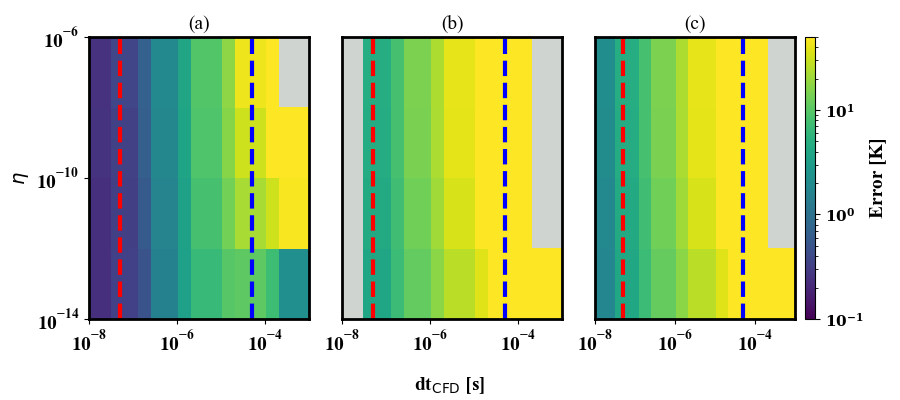}
    \includegraphics[width=0.47\textwidth]{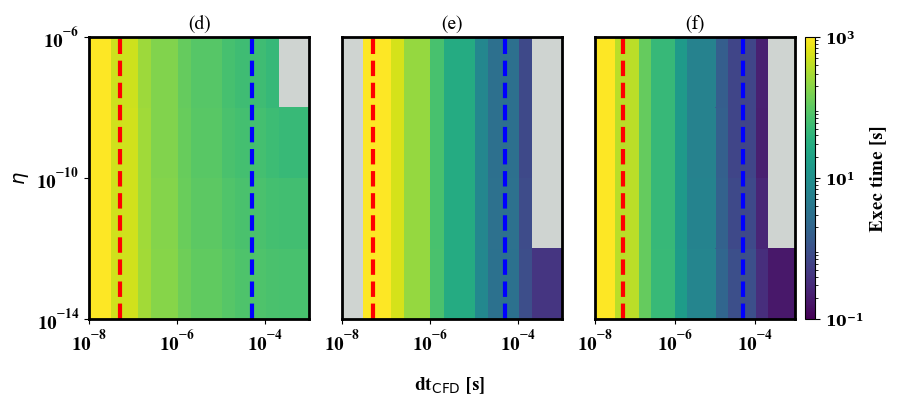}
    \caption{Using a direct linear solver with either a numerically computed or an analytically computed Jacobian is more robust than a Newton-Krylov (NK) approach for $dt_{CFD}$ in the PeleLMeX regime (near the blue dashed line). I.e., it allows the integrator to meet the error requirements with a timestep that is larger and therefore the overall integration is faster than with NK and GMRES.
    In the PeleC regime (near the red dashed line) NK with GMRES is effective due to the small time scales.
    The top row shows the averaged mean squared error for the 0D reactor test problem when using solution approaches 1B, 2B, and 3B from Table~\ref{tab:sol_app} respectively.
    The bottom row shows the execution time.
    In this case, gray indicates the run timed out, or that the error was so large that the reactants did not ignite.}
    \label{fig:compare_solves}
\end{figure}

We achieved further performance improvements in PeleLMeX by using direct linear solvers within a modified Newton approach.
The smaller time scales of PeleC result in SUNDIALS selecting its time steps based on accuracy considerations (rather than solver robustness), and the resulting linear systems tend to be diagonally dominant. Thus, PeleC problems tend to allow the effective use of GMRES without a preconditioner to aid convergence. However, the larger time scales of PeleLMeX result in linear systems that are more difficult to solve with unpreconditioned Krylov iterations. In the PeleLMeX case, the choice of the linear solver used within the Newton iterations can be critical to the robustness of the overall time integration.



When unpreconditioned GMRES is employed as the linear solver, the difficult linear systems that arise when integrating over PeleLMeX timescales lead to failed SUNDIALS ODE integrations.  They either time out (according to the integrator settings supplied), or they fail to ignite (Figure~\ref{fig:err_sol}).
Specifically, we observe such failures when $dt_{CFD}$ is near $10^{-5}$ or larger.
The failures are due to at least two compounding effects.
First, if the individual components are not scaled properly, as discussed above, we may need to converge the iterated systems to lower tolerances to ensure stable and accurate system evolution. Note that although we can modify the settings to account for this, we also observe that as the integration intervals ($dt_{\textit{CFD}}$) increase, the chemical state can vary considerably and complicate choosing an effective scaling matrix that remains suitable over the entire interval.
Second, as the size SUNDIALS-selected time steps increase, the underlying linear systems lose diagonal dominance.
In these cases, without effective preconditioning, the Krylov solvers will require an increasing number of iterations to converge the linear system.
Eventually, the linear solver will hit a user-specified maximum allowable iterations and  trigger SUNDIALS to cut its internal time step size and retry the step.
Thus, for any reasonable configuration, the integration eventually times out. Although preconditioning such systems is a typical strategy to reduce Krylov iterations, development of an effective, suitably general, memory efficient preconditioner for combustion chemistry networks
that is also suitable for application on GPUs remains an area of ongoing research.

Using the zero-dimensional reactor problem problem, we compare the performance of a dense direct linear solver (based on LU factorization with pivoting) paired with a numerical Jacobian and with an analytical Jacobian.
The numerical Jacobian is formed using difference quotients as described in \citep[\S2.1]{hindmarsh2005sundials}.
This procedure requires an additional call to the ODE right-hand-side per column of the Jacobian and a well-suited perturbation in order to form an accurate numerical partial derivative for each entry.
The analytical Jacobian is computed using model-specific code that is generated offline.

Relative to the Newton-Krylov option, modified Newton with a dense direct linear solver paired with a numerical or analytical Jacobian in the modified Newton solver strategy is at least one order of magnitude faster for $dt_{CFD}$ in the PeleLMeX regime.
This is because these solver strategies allow the integrator to take larger timesteps. In the PeleC regime, the Newton-Krylov option is both more accurate and faster than the modified Newton strategies.
\addressed{MH}{}{I think we want to be explicit here instead of saying "robust": The Newton solver is orders of magnitude faster for LMeX type of steps, albeit a bit more inaccurate. At C type of steps, the Newton-Krylov method is more accurate and faster}
For the cases tested, the analytical Jacobian is consistently faster than the numerical Jacobian (Figure~\ref{fig:compare_solves}).

\subsubsection{Explicit Methods}

\begin{figure}
    \centering
    \includegraphics[width=0.45\textwidth]{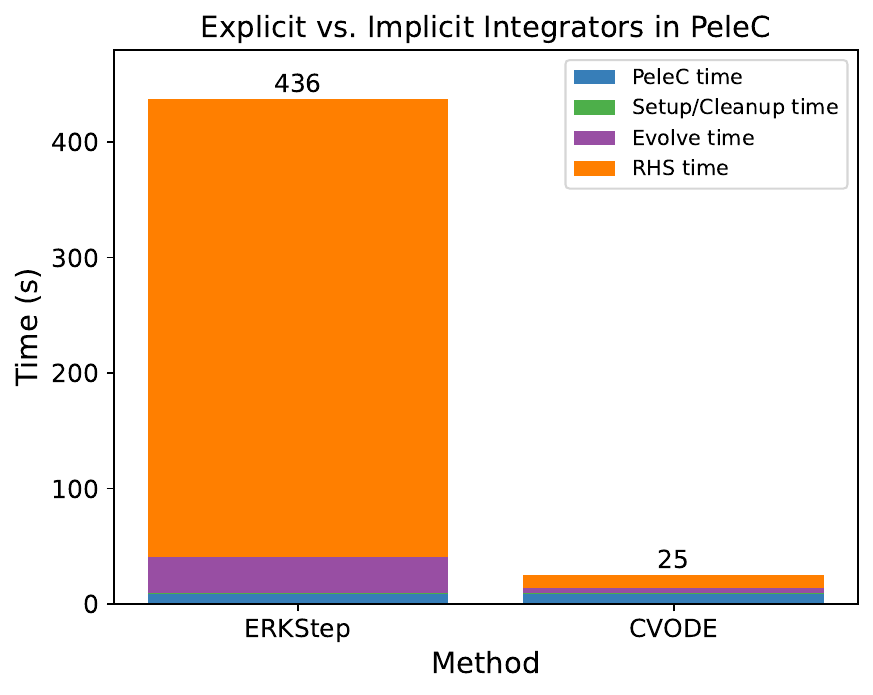}
    \caption{Comparison timings for an explicit method from ARKODE and the implicit methods from CVODE on a pre-mixed flame test case in PeleC. While the explicit method is less expensive per step, the greater stability of the implicit method enables it to take far fewer steps leading to a significantly faster runtime.}
    \label{fig:pelec_exp_v_imp}
\end{figure}

Depending on the stiffness of the batched systems, the explicit time integration methods provided by ARKODE can be a potential alternative to implicit approaches in CVODE. These methods provide similar adaptivity in internal time step size, but employ fixed order schemes. Importantly, they do not require algebraic solvers, removing much of the complexities discussed above. However, multistage explicit methods may require more work per step than the linear multistep methods in CVODE, depending on the performance of the algebraic solvers selected for the specific problems being addressed.

ARKODE and CVODE are built on same shared SUNDIALS infrastructure with nearly identical user interfaces, and applications can easily interface to both codes to directly compare the performance of different classes of methods on specific problems. Figure~\ref{fig:pelec_exp_v_imp} shows a breakdown of the timings from a simulation with a fourth-order explicit method and solution approach 1B with CVODE on a pre-mixed flame test in PeleC with a 53 species Dodecane chemistry mechanism. The domain has size 0.4 × 0.4 × 1.6 with a base grid of size 32 × 32 × 128 cells and one level of mesh refinement in the center of the domain around the flame, with a refinement ratio of two. Ten coarse time (hydrodynamic) time steps are taken starting from an interpolated initial condition. While the explicit method requires slightly less work on average per step than the implicit method, $5.2$ RHS evaluations per step compared to $6.0$, it requires far more steps on average, $196.3$ compared to $4.9$, as the implicit method can step over fast dynamics the explicit method tracks. Thus, for this problem and combustion model the implicit methods in CVODE significantly outperform the explicit approach.


\subsubsection{Batched Linear Solvers}

\addressed{CSW}{}{It would be good to cite the Ginkgo sundials paper somewhere in this section.  Like, "Another option to MAGMA is to use a sparse iterative approach for the linear systems through Gingko.  We investigated this approach, and while it gave some advantages, ultimately it was not as fast as the dense solvers in MAGMA for the specific reaction systems we were most interested in for the Exascale Computing Project goals.  Further information on that work can be found here."  Or something along lines like that. CJB: Ginkgo was actually was faster for the interesting reaction mechanisms (according to there most recent testing), but I added a sentence along these lines.}

\begin{figure}[!htb]
\centering
\includegraphics[width=0.45\textwidth]{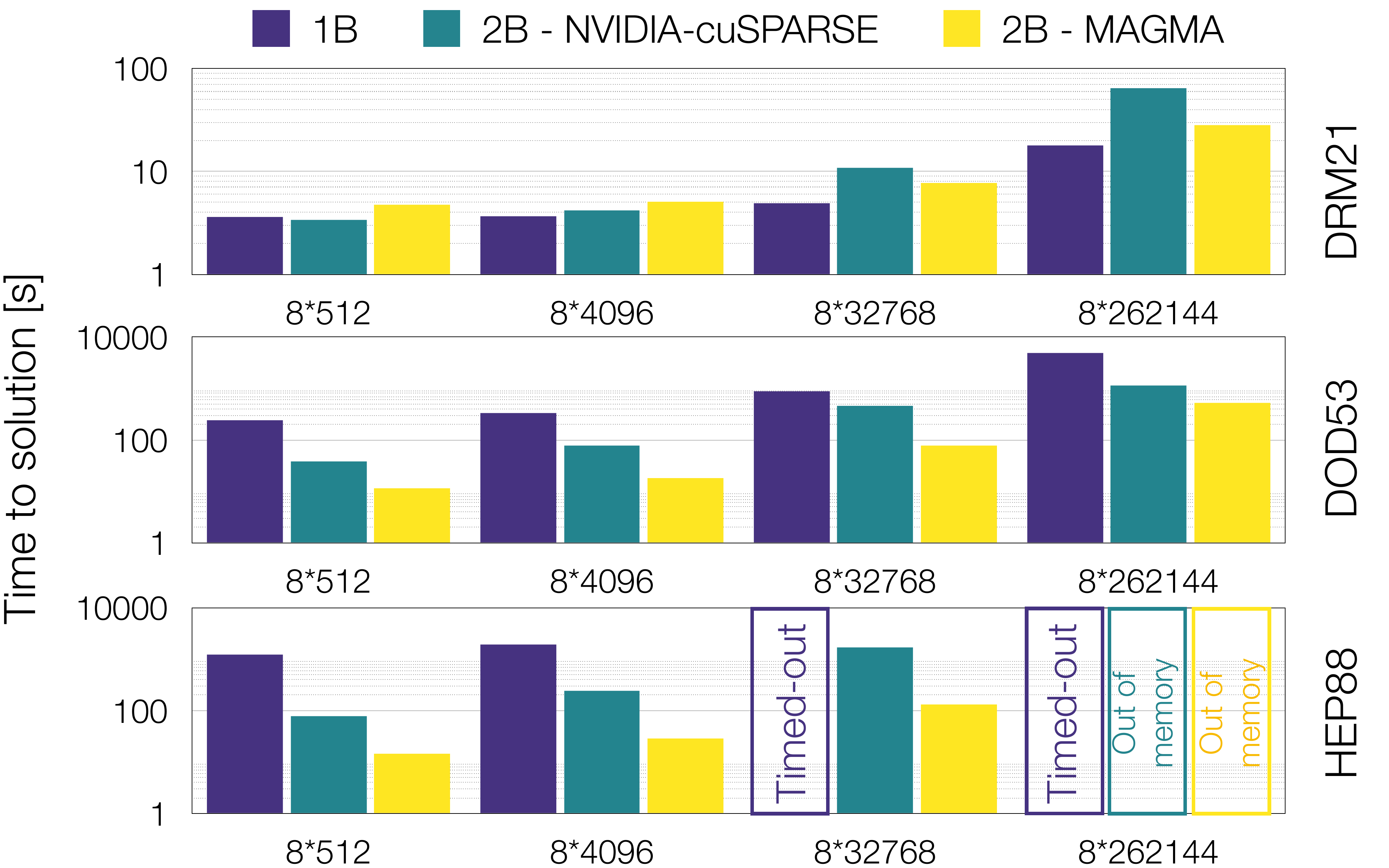}
\caption{Batched direct linear solvers that exploit the independence of the reaction ODE in each mesh cell are faster than using GMRES for medium (DOD53) and large (HEP88) Pele chemical mechanisms. GMRES is still more efficient in most cases with a small mechanism (DRM21). Which solver is best also depends on the batch size (512, 4096, 32768, or 262144).}
\label{fig:pele_batched_linsolve}
\end{figure}

As demonstrated above, using a direct linear solver in CVODE can be beneficial to the ODE integrator robustness. On GPUs, direct solvers can prove inefficient and extremely memory consuming for large linear systems. It is possible to exploit the linear system block-diagonal pattern presented in Figure~\ref{fig:pele_ordering_PP}(a) by using a batched linear solver, where the individual linear systems corresponding to each cell are considered independently. To evaluate the effect of the linear solver choice, we compare approach 1B with two variants of approach 2B (see Table \ref{tab:sol_app}).  In the first variant, the batched direct solve is performed with the NVIDIA cuSPARSE library \citep{cusparse}, while the MAGMA library  \citep{magma} is employed in the second case.
Another option is to use a sparse iterative approach for the linear systems through the Ginkgo library \citep{Anzt_Ginkgo_A_Modern_2022}, but this is primarily useful for reaction mechanisms larger than the ones considered here \citep{aggarwal2021batched}.
The test case consists of a rectangular computational domain filled with a uniform mixture of air and hydrocarbon, and an increasing temperature profile in one direction such that part of the domain will experience auto-ignition during the course of the simulation. Three chemical mechanisms of increasing complexity, in terms of number of species and reactions, are considered.
The domain is divided into 8 grid patches that range in size from 8$^3$ to 64$^3$ to evaluate the effect of the batch size (i.e., number of grid cells in a grid patch). The time to solution is presented in Figure~\ref{fig:pele_batched_linsolve}, comparing the three linear solver approaches for each mechanism and batch size. No single method is found to outperform all the others across the full set of test cases. Rather, approach 1B is consistently more efficient on this problem for the smaller chemical mechanism, whereas for larger mechanisms, both variants of 2B are significantly faster than 1B. Note that the 2B-MAGMA case is 1.5 to 10 times faster than its cuSPARSE counterpart.  We also note that memory requirements of the larger models can limit the range of grid patch sizes available.  Finally, since the batched MAGMA solvers are considerably more efficient than the cuSPARSE counterparts, we will exclusively use the MAGMA solvers for the remaining examples in this paper.

\subsection{Improving performance in Nyx}
\label{ss:improving_nyx}

In this section we discuss three algorithm and software changes implemented to improve the performance of the ODE integrations in Nyx.

\subsubsection{Increasing Concurrency with OpenMP and Streams}
\begin{figure}
\centering
\includegraphics[width=.5\textwidth,keepaspectratio,trim={8mm 14mm 8mm 14mm},clip]{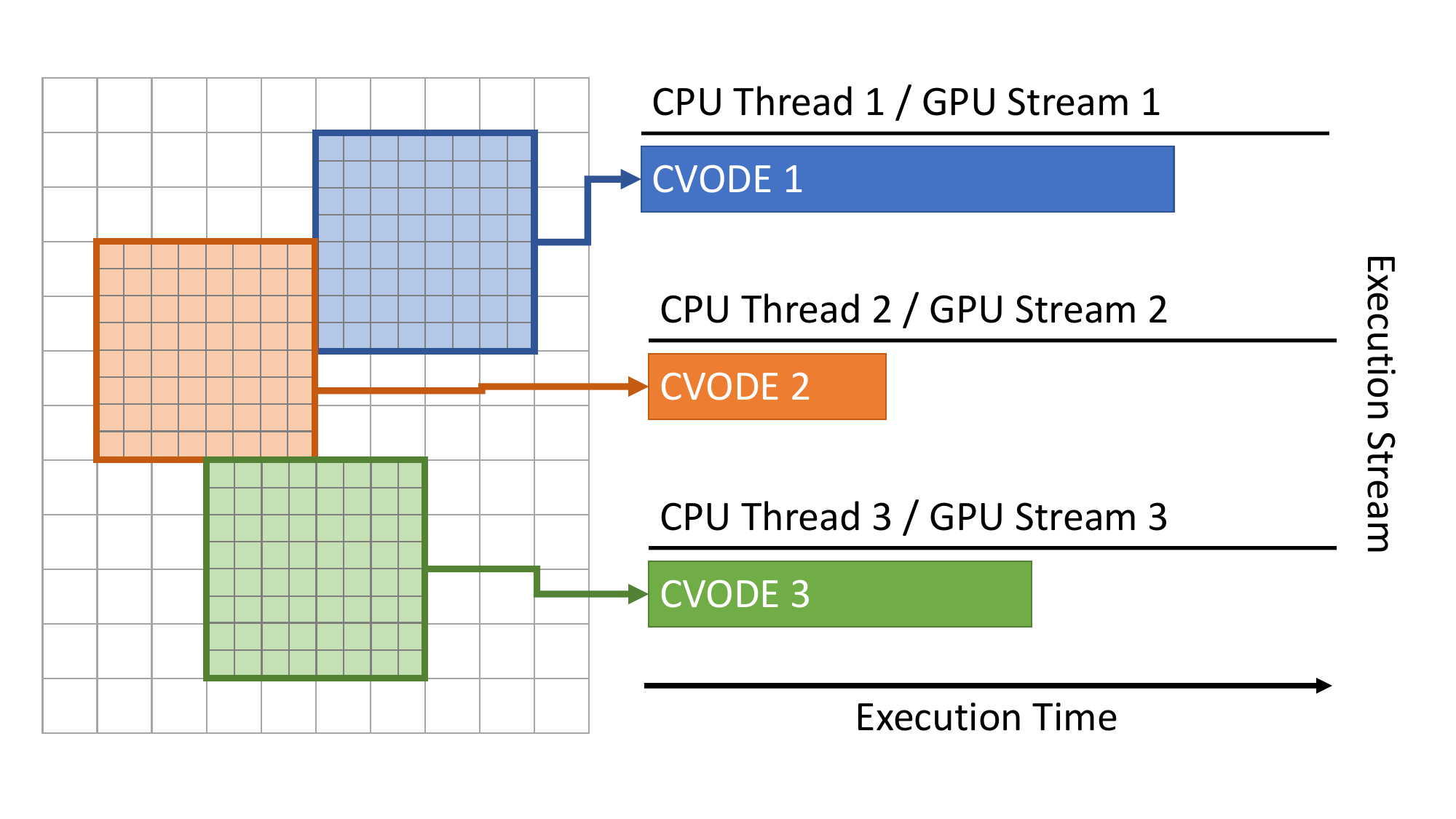}
\caption{Diagram illustrating the additional concurrency leveraged when utilizing multiple integrator (CVODE) instances to simultaneously evolve different boxes of cells, each associated with a unique CPU thread and GPU stream.}
\label{fig:cpu_gpu_diagram}
\end{figure}

Since only a single scalar ODE is evolved per spatial grid cell in Nyx, increasing concurrency in the ODE evolution is essential to achieving good performance on GPU systems. We achieve this concurrency by batching cells together. However, the AMR configuration determines the maximum number of cells in a box and limits the maximum number of cells available to the integrator.
Thus, to enable greater parallelism, we create multiple SUNDIALS instances, each one concurrently advancing the box of cells assigned to that CPU thread and GPU stream as illustrated in Figure~\ref{fig:cpu_gpu_diagram}. Once an integrator instance completes the evolution of a box of cells, it is given a new box to advance until all the work is completed.

\begin{figure}
\centering
\includegraphics[width=.5\textwidth,height=\textheight,keepaspectratio,trim={3mm 4mm 3mm 0mm},clip]{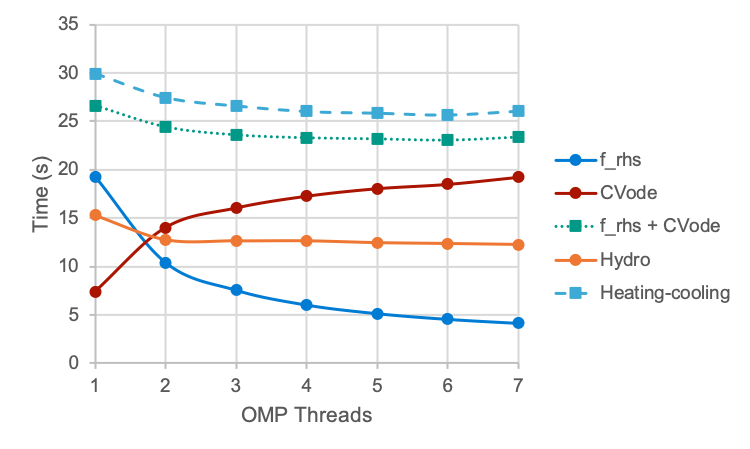}
\caption{Utilizing concurrent CVODE instances (one per OpenMP thread) leads to a 14\% improvement in performance with 6 threads compared to a single instance starting from a physically interesting time (redshift $z=7$) in a Nyx Lyman-$\alpha$ example with $256^3$ cells taking 50 steps and run on NVIDIA V100 hardware.}
\label{fig:nyx_timings_threading}
\end{figure}

Figure \ref{fig:nyx_timings_threading} shows the timings from different parts of Nyx with increasing numbers of OpenMP threads, i.e., a growing number of CVODE instances running currently, starting from a physically interesting time (redshift $z=7$) in a Nyx Lyman-$\alpha$ example with $256^3$ cells taking 50 steps and run on NVIDIA V100 hardware. Overall, the time spent in the heating-cooling computation decreases with diminishing returns as the number of threads (and CVODE instances) increases. The greater concurrency leads to a 14\% reduction in run time with 6 threads compared to a single execution stream.


\subsubsection{Kernel Fusion}
\begin{figure}
\centering
\includegraphics[width=.49\textwidth,height=.75\textheight,keepaspectratio, trim={3mm 0mm 3mm 0mm},clip]{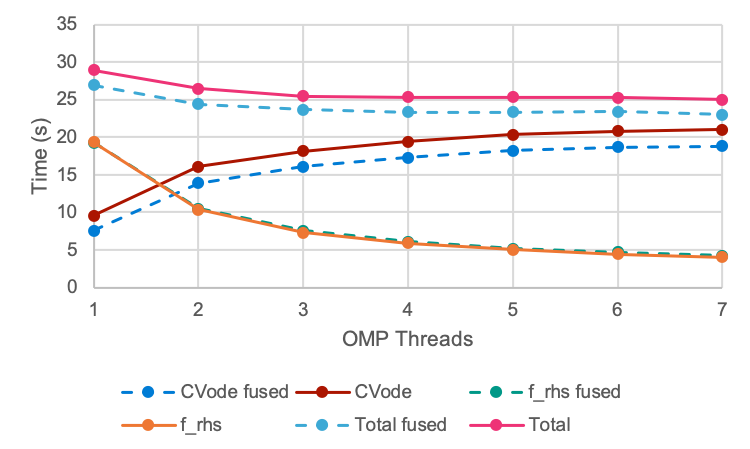}
\caption{Utilizing fused kernels in CVODE leads to a 7\% to 8\% improvement in the total run time for a Nyx Lyman-$\alpha$ test with $256^3$ cells taking 50 steps from redshift $z=7$ on NVIDIA V100 hardware.
}
\label{fig:nyx_timings_fusing}
\end{figure}

Several computations within CVODE require performing multiple vector operations in succession to compute a needed result.
For example, within the diagonal linear solver utilized by Nyx, setting up and inverting the linear system within the Newton solver requires ten vector operations in succession.
Each of these operations necessitates a kernel launch and pass over all of the vector data.
We enhanced CVODE with fused versions of many common sets of operations to reduce the kernel launch overhead and number of memory accesses.
Figure \ref{fig:nyx_timings_fusing} shows the timings with and without fused integrator kernels for various numbers of concurrent CVODE instances (one per OpenMP thread) on the same Lyman-$\alpha$ example from the prior section.
Enabling fused kernels leads to a consistent improvement of approximately 7\% to 8\% in the total run time regardless of the number of CVODE instances utilized.



\subsubsection{Tiling}
An important consideration when optimizing for different architectures is tuning the problem decomposition across MPI ranks, OpenMP threads, and GPU streams with the goals of minimizing runtime and running the largest simulation possible given limited computational resources. Depending on the box size, the amount of GPU scratch memory space available to computational kernels can be the limiting factor on the maximum problem size per node.
To address this issue, we leverage the tiling capability in AMReX \citep{zhang2021amrex} to transform loops over each grid patch into partitioned iterators that loop over subregions that tile the patch. Because each subregion is updated independently, multiple CVODE instances can work on the patch simultaneously in their own GPU stream with no race conditions.
These partitioned iterators enable better utilization of the available high bandwidth memory and increase concurrency.

\begin{figure}
\centering
\includegraphics[width=.49\textwidth,height=.75\textheight,keepaspectratio,trim={0mm 0mm 0mm 0mm},clip]{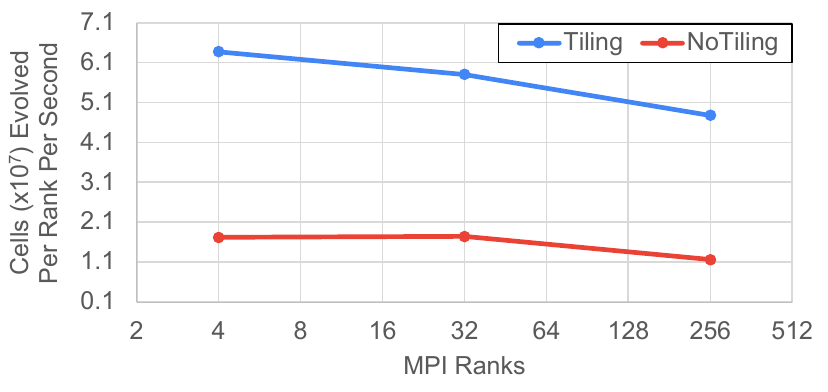}
\caption{
Performance with and without tiling from a Nyx Lyman-$\alpha$ simulation with ${512}^3$ to ${2048}^3$ cells and particles on NVIDIA A100 hardware (Perlmutter) taking 100 steps from redshift $z=2.2$ using 8 OpenMP threads and CUDA streams with $2*{256}^3$ cells per rank.
In this setup the tiling can increase throughput by approximately 3$\times$ to 4$\times$ by enabling additional concurrency.
}
\label{fig:nyx_timings_tiling_pm}
\end{figure}

Figure~\ref{fig:nyx_timings_tiling_pm} shows the number of cells advanced per rank per second, which divides the total number of cells by the number of MPI ranks and by the average time to complete a $dt_{CFD}$ time step (higher values indicate better performance), in a Lyman-$\alpha$ simulation. The simulations were run with and without tiling using ${512}^3$ to ${2048}^3$ cells and particles ($2*{256}^3$ cells per rank) on NVIDIA A100 hardware taking 100 steps from redshift $z=2.2$ with 8 OpenMP threads and CUDA streams and fused operations. Tiling enables significantly faster results with the same resources due greater concurrency producing a 3$\times$ to 4$\times$ increase in the number of cells advanced per rank per second.

\begin{figure}
\centering
\includegraphics[width=.49\textwidth,height=.75\textheight,keepaspectratio,trim={0mm 0mm 0mm 0mm},clip]{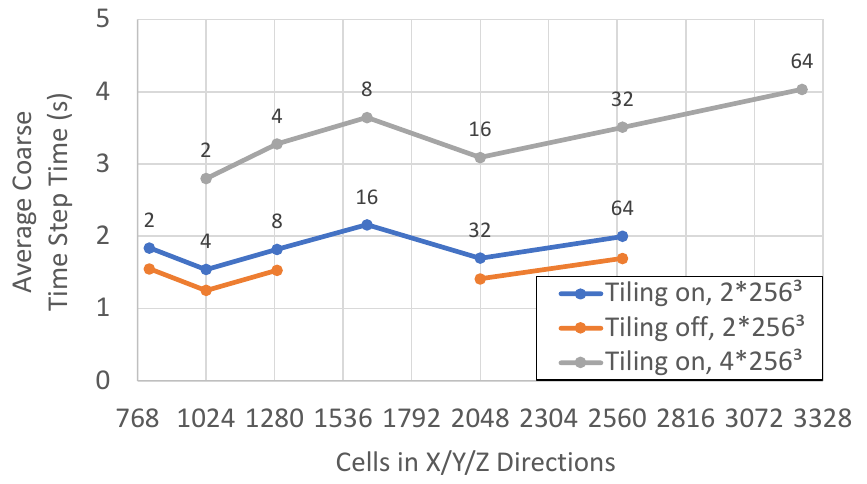}
\caption{Average time to complete a coarse step ($dt_{CFD}$) with $2 * 256^3$ or $4 * 256^3$ cells per node with and without tiling on a Nyx-Lyman-$\alpha$ simulation taking 120 steps from redshift $z=200$ on 2 to 64 nodes (annotations) with AMD MI250X hardware. Runs utilize a Nyx option limiting concurrency to make more scratch memory available to hydrodynamic kernels. This increases the run time when sufficient memory is available ($2 * 256^3$ cells per node) but enables running larger problems ($4 * 256^3$ cells per node) than possible on the same resources without tiling. Note that the problem setup on 16 nodes without tiling did not complete successfully.
}
\label{fig:nyx_timings_tiling}
\end{figure}

Alternatively, tiling can be used to run larger problem sizes with a given set of resources. Figure~\ref{fig:nyx_timings_tiling} shows the average time to complete a $dt_{CFD}$ time step with two different problem sizes per node, with and without tiling, in a Nyx Lyman-$\alpha$ simulation with ${768}^3$ to ${3328}^3$ cells and particles on AMD MI250X hardware taking 120 steps starting from redshift $z=200$. The simulations use 8 OpenMP threads and have fused kernels enabled. Additionally, a Nyx option is enabled to reduce the concurrency and make more scratch memory available in hydrodynamic computations by limiting the code to working only on one tile at a time. When the problem is decomposed with $2 * 256^3$ cells per rank, tiling increases the average time per step by approximately 20\%. However, because tiling reduces the required resources per kernel launch, it enables running with a larger problem size per rank than possible without tiling, $4 * 256^3$ cells.

\subsection{Improving performance: common themes}
\label{ss:improving_common}

\subsubsection{Efficient Memory Usage}
\begin{figure}
\centering
\includegraphics[width=.45\textwidth,height=.75\textheight,keepaspectratio]{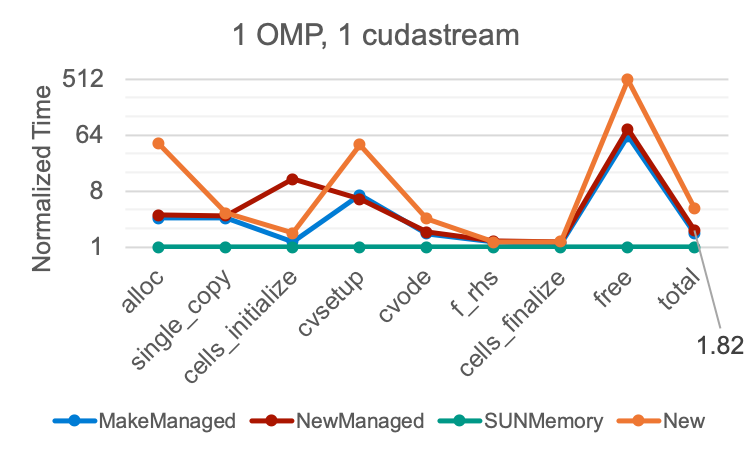}
\caption{Performance of different profiling regions for various allocation strategies in a Nyx Lyman-$\alpha$ $256^3$ simulation on NVIDIA V100 hardware taking 50 steps from redshift $z=7$. Times are normalized against the setup using the SUNMemory API which performed allocations through an interface to the AMReX memory pools. With the ``New'' approach SUNDIALS performs host and device allocations. With ``NewManaged'', SUNDIALS performs managed memory allocations rather than host and device memory allocations. With ``MakeManaged'', the initial state memory is provided by Nyx but all other allocations are performed by SUNDIALS. The SUNMemory approach removes artificial competition for memory and reduces the number of actual memory allocations/de-allocations leading to an approximately $2\times$ speedup in the total Nyx run time compared to the ``Managed'' approaches.
}
\label{fig:nyx_timings_allocation}
\end{figure}

One issue common to both Pele and Nyx is efficient memory usage and data exchange between the application code and SUNDIALS. Any application that utilizes SUNDIALS exchanges data with it through a vector data structure that contains the initial condition of the ODE.
SUNDIALS takes the initial condition vector and ``clones'' it, i.e. creates a new instance of a vector with the same properties but a new and distinct data array, to create all of the internal vectors it needs to integrate the problem.
Since Pele and Nyx implicitly and explicitly use large pools of CPU, GPU, and other memory allocated by the AMReX framework at initialization, an artificial competition between SUNDIALS and AMReX for memory resources arises.
The competition is artificial since Pele and Nyx do not utilize all of the memory available, but it is considered allocated to the AMReX memory pool by the operating system and/or device driver.
To solve this problem, the ``SUNMemory'' API was added to SUNDIALS to support application memory pools \citep{balos2021enabling}, and a corresponding interface to it was added within AMReX. This interfacing ensures that any ``clones'' of the initial condition vector in SUNDIALS draw from the AMReX memory pool.
Furthermore, the data is kept resident in the GPU device memory throughout.
This residence is critical to runtime performance since moving data from CPU memory to GPU device memory is costly.
In addition to solving an artificial competition for memory, the use of the SUNMemory API improves runtime performance by eliminating extraneous calls to allocate memory every time Pele and Nyx initialize SUNDIALS (which occurs at every fluid time step).
A total runtime performance improvement of approximately $2\times$ in Nyx was observed (Figure \ref{fig:nyx_timings_allocation}) when using the SUNMemory API compared to other allocation approaches.


\subsubsection{Fast and Flexible Reductions}

\begin{figure}[!htb]
\centering
\includegraphics[width=0.45\textwidth]{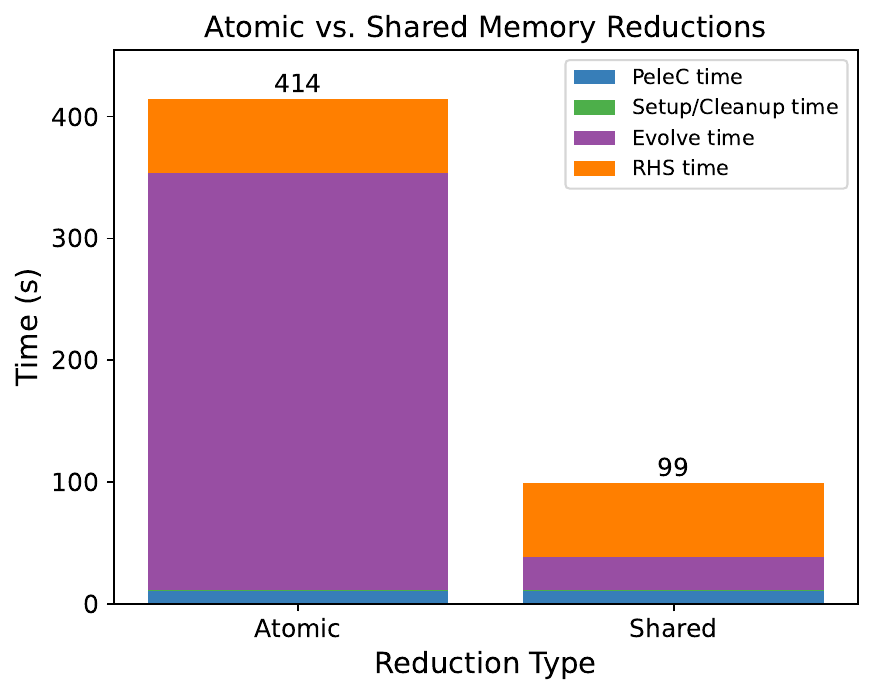}
\caption{PeleC timings with CVODE on a pre-mixed flame (PMF) test problem on AMD MI100 GPUs using approach 1a with vector reduction operations in SUNDIALS implemented using either atomic operations (left) or shared memory (right). As the MI100 does not have hardware support for double precision atomics, the shared memory implementation significantly reduces the overall runtime.}
\label{fig:pelec_reductions}
\end{figure}

Reduction operations are common within the SUNDIALS integrators and arise in the computation of vector norms, which are used for error estimates in time step selection and in the stopping criteria of iterative algebraic solvers, and dot products, which occur frequently within the iterative linear solvers.
There are two approaches to implementing these reduction operations for GPU devices supported in SUNDIALS that Pele and Nyx both leverage.
The first approach utilizes fast ``shared'' GPU memory to accumulate values while the second approach utilizes atomic operations.
Depending on the specific GPU architecture targeted, one of these approaches may be more efficient than the other. Specifically, when atomic operations are not backed by hardware support, it becomes critical to use the ``shared'' memory approach.
We investigated the performance difference by comparing the timings for two variations of solution approach 1A in a PeleC pre-mixed flame test problem, one using atomic operations and the other using the ``shared'' memory approach. With this solution approach, where GMRES is the linear solver, efficient reduction operations are critical as an increasing number of dot products are computed in the modified Gram-Schmidt algorithm each linear solver iteration within each nonlinear solver iteration.
The runs utilized AMD MI100 GPUs, which do not have hardware support for atomic operations in double precision, and show that the ``shared'' memory approach is approximately four times faster than using atomic operations on this hardware (Figure~\ref{fig:pelec_reductions}).
With hardware that does support the needed atomic operations, we have observed the atomic reduction approach can be faster than the shared memory approach.
As such, SUNDIALS allows for codes to switch between the two approaches.


\section{Results and Discussion}\label{s:results}


We now demonstrate the performance of Pele, Nyx, and SUNDIALS on problems that exhibit sufficient complexity to require the use of many of the advanced SUNDIALS features and performance improvements we have discussed.
These problems are run on the Summit petascale and Frontier exascale supercomputers at the Oak Ridge Leadership Computing Facility (OLCF), and the Perlmutter supercomputer at the National Energy Research Scientific Computing Center (NERSC).
Summit is a 200 petaflop (peak) IBM machine comprised of 4608 nodes each with two 22-core IBM Power9 CPUs and six NVIDIA V100 GPUs. Each node contains 512 GB of DDR4 memory and 96GB of High Bandwidth Memory (HBM2) \citep{SummitUserGuide}.
Frontier is a 1.6 exaflop (peak) Cray EX machine with 9408 nodes. Each node of Frontier has a 64-core AMD ``Optimized 3rd Gen EPYC'' CPU connected to 512GB of DDR4 memory and 4x AMD MI250x GPUs with 128GB of HBM2E memory \citep{FrontierUserGuide}.
Perlmutter is a 71 petaflop (peak) Cray EX machine with 3072 CPU-only and 1792 GPU-accelerated nodes. In this work we use the GPU-accelerated nodes which consist of a 64-core AMD EPYC 7763 CPU paired with 256GB of DDR4 memory and 4x NVIDIA A100 GPUs with up to 80GB of HBM2E memory \citep{NerscUserGuide}.

\subsection{Pele}
Our Pele demonstration utilizes PeleLMeX for a multipulse injection scenario, typical of modern combustion systems, where a short pulse of a high-reactivity (n-dodecane) fuel jet enters a high pressure reaction chamber that is prefilled with a low-reactivity air-fuel mixture (methane-air). As the fuel jet develops, it fills the domain with a highly turbulent combustible gas mixture. After a short dwell time, a second pulse injects another stream of high-reactivity diesel fuel (n-dodecane) into the chamber. The second pulse interacts with the first pulse that is already undergoing a low temperature autoignition process that is followed by a high-temperature flame propagation stage.   Autoignition is a local process dominated by chemistry only, whereas flame propagation involves a detailed balance between advection, diffusion, and chemical processes. In this system, the combustion processes span timescales that are far less than to far greater than the timescales of turbulent transport near the fuel jets. Thus, this system requires the use of the stiff integrators provided by SUNDIALS.  The chemistry model contains 53 species that vary in peak mass fractions from ${\cal O}$(1e-10) to ${\cal O}$(1), and temperature varies over space and time from 450K to over 2400K. Also, much of the domain is filled either with an unreacted fuel-air mix, or with hot products at relatively uniform temperatures. After ignition, the solution contains highly localized pockets with significantly heterogeneous distributions of species, temperature, and chemical activity.

\begin{figure}
    \centering
    \begin{subfigure}{0.24\textwidth}
        \centering
        \includegraphics[width=\textwidth]{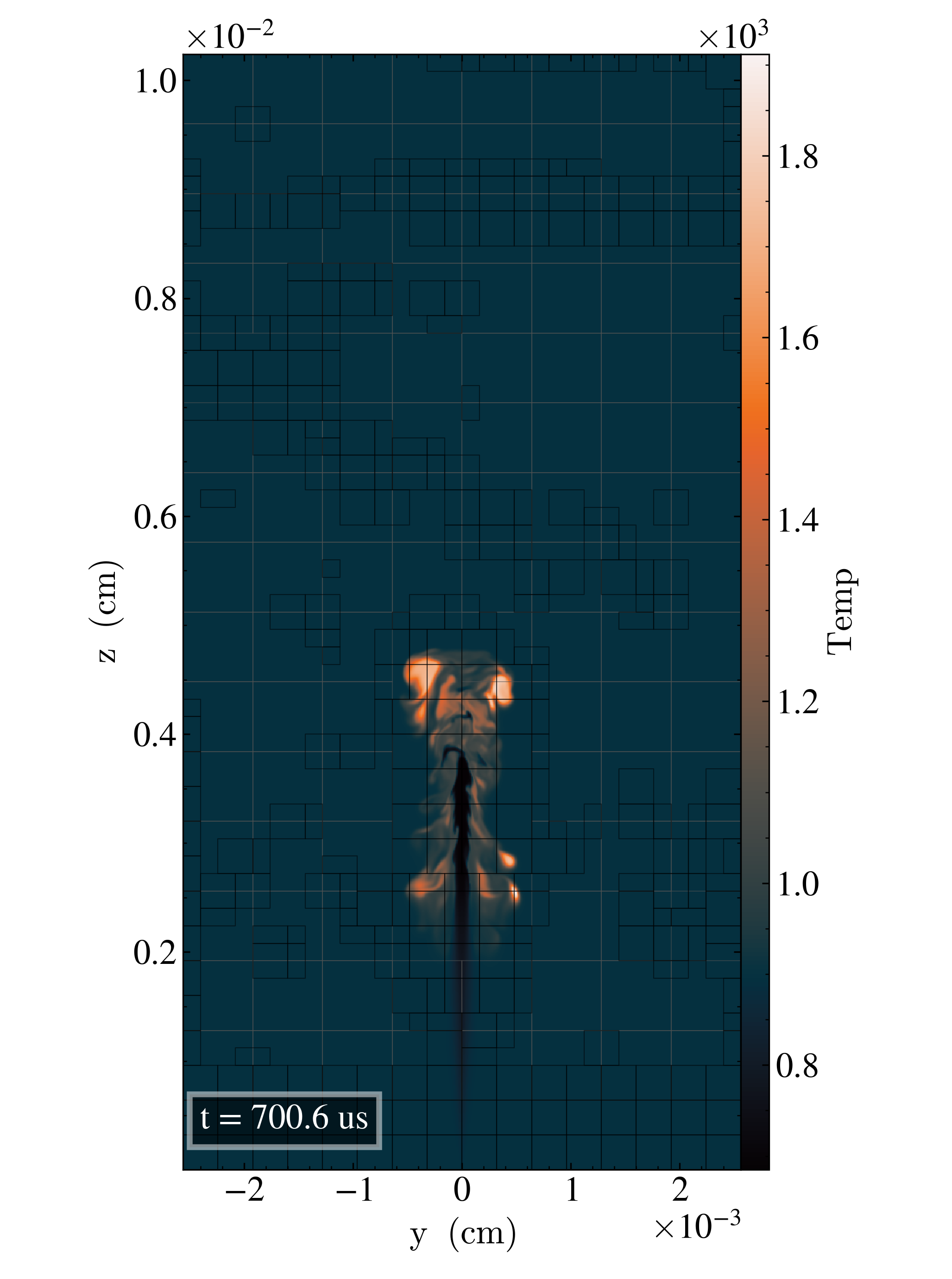}
        \caption{Start of ignition}
        \label{subfig:pelelmex_start}
    \end{subfigure}
    \begin{subfigure}{0.24\textwidth}
        \centering
        \includegraphics[width=\textwidth]{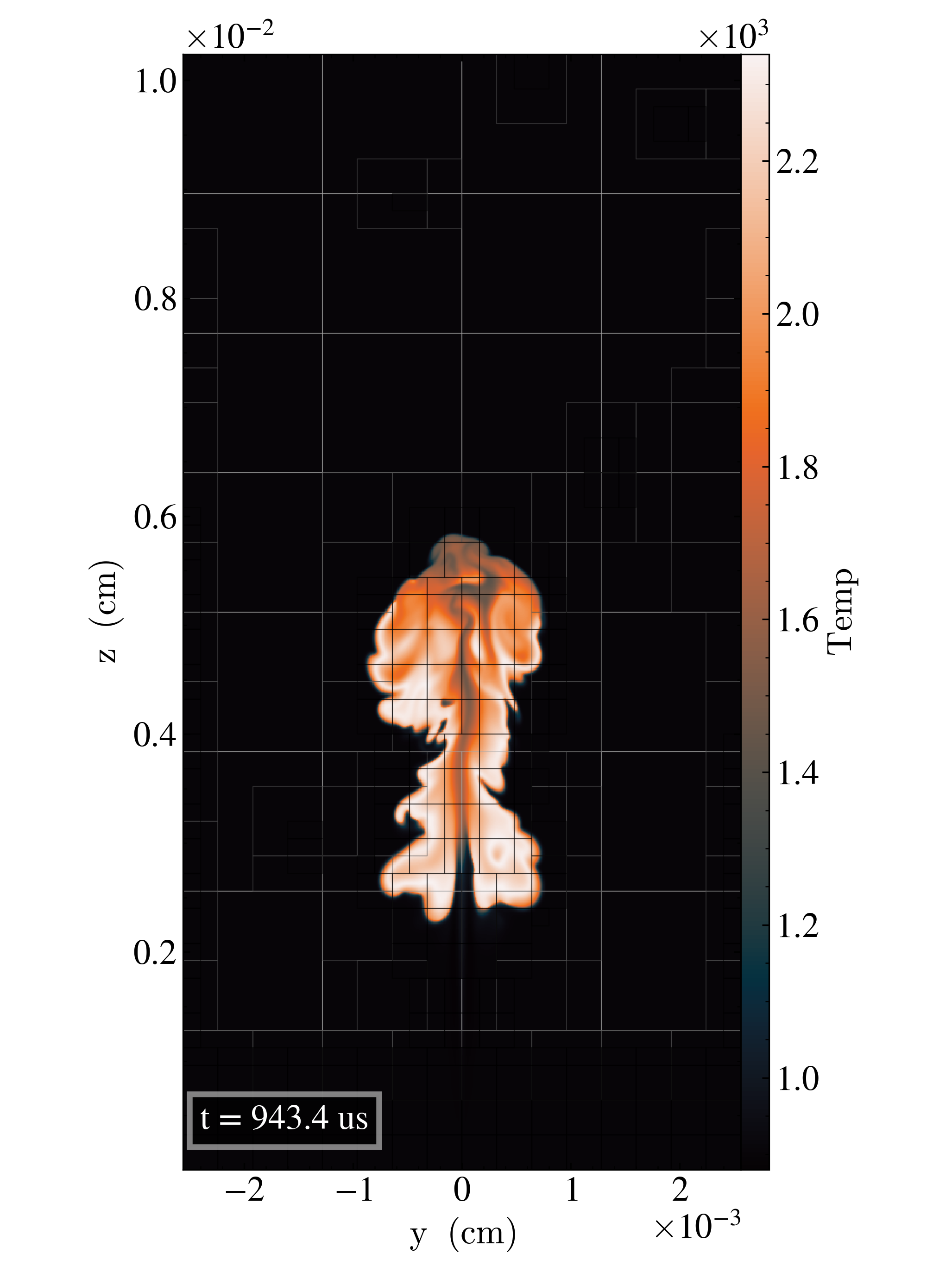}
        \caption{Fully-developed ignition}
        \label{subfig:pelelmex_develop}
    \end{subfigure}
    \caption{PeleLMeX pulsed jet injection case showing a slice of the temperature field near the start and development of ignition.}
    \label{fig:pelelmex_ignition}
\end{figure}

Figure~\ref{fig:pelelmex_ignition} shows the development of the first fuel injection pulse of n-dodecane into the methane-air premixture. Figure~\ref{subfig:pelelmex_start} shows the start of the auto-ignition process with the high-reactivity diesel fuel developing ignition kernels near the shear mixing layer with the high temperature background methane-air premixture. Figure~\ref{subfig:pelelmex_develop} shows the development of the higher temperature flame propagation with the start of the second injection that will soon penetrate into the hot combustion products that are produced by the first injection. At this stage of the combustion process, the evaluation of the chemical reaction mechanism is a significant portion of the overall CFD run-time thereby motivating the need for efficient algorithms and scalability provided by SUNDIALS.


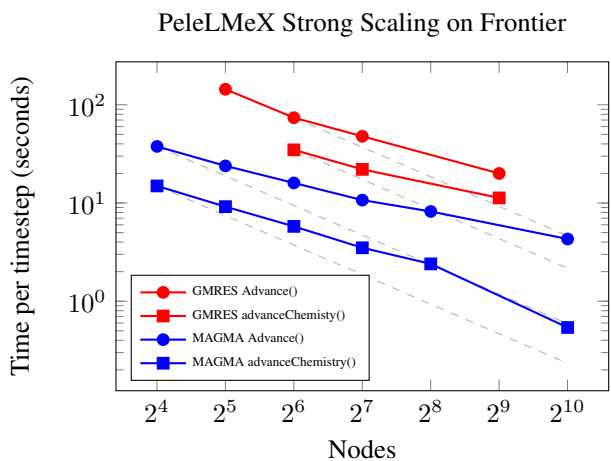
\begin{figure}
\centering
\begin{tikzpicture}
\begin{loglogaxis}[
  title={PeleLMeX Strong Scaling on Frontier},
  xlabel={Nodes},
  ylabel={Time per timestep (seconds)},
  log basis x={2},
  legend pos=south west,
  legend columns=1,
  legend style={font=\tiny},
  legend cell align={left},
  width=0.95\linewidth,
  height=0.7\linewidth,
  every axis plot/.append style={thick},
  /pgf/number format/.cd,
        use comma,
        1000 sep={},
]

\addplot[thin, black!30, mark=none, dashed, forget plot] coordinates {
  (16, 37.7)
  (1024, 0.59)
};

\addplot[thin, black!30, mark=none, dashed, forget plot] coordinates {
  (16, 14.92)
  (1024, 0.233)
};

\addplot[thin, black!30, mark=none, dashed, forget plot] coordinates {
  (64, 73.95)
  (1024, 4.62)
};

\addplot[thin, black!30, mark=none, dashed, forget plot] coordinates {
  (64, 34.78)
  (1024, 2.17)
};

\addplot[red, mark=*, mark options={solid}] coordinates {
  (32, 144)
  (64, 73.95)
  (128, 47.8)
  (512, 20.0)
};

\addplot[red, mark=square*, mark options={solid}] coordinates {
  (64, 34.78)
  (128, 22.0)
  (512, 11.27)
};

\addplot[blue, mark=*, mark options={solid}] coordinates {
  (16, 37.7)
  (32, 23.88)
  (64, 16.0)
  (128, 10.7)
  (256, 8.2)
  (1024, 4.3)
};

\addplot[blue, mark=square*, mark options={solid}] coordinates {
  (16, 14.92)
  (32, 9.18)
  (64, 5.78)
  (128, 3.5)
  (256, 2.4)
  (1024, 0.54)
};

\legend{GMRES Advance(), GMRES advanceChemisty(), MAGMA Advance(), MAGMA advanceChemistry(),}
\end{loglogaxis}
\end{tikzpicture}
\caption{Strong scaling of the n-dodecane skeletal reaction mechanism with PeleLMeX using GMRES (red lines) or MAGMA with an analytic Jacobian (blue lines) on 128 to 1024 Frontier nodes with 78M cells on three AMR levels. The gray dashed lines indicate the ideal scaling from the initial point.}
\label{fig:pelelmex_strong_scaling}
\end{figure}

The 53 species n-dodecane chemical reaction was evaluated using both approach 1B and 2B (and leveraging MAGMA, as discussed above). To evaluate performance of each of the solvers, we performed a strong scaling study using the PeleLMeX simulation shown in Figure~\ref{fig:pelelmex_ignition}.
This case has 78 million grid cells and uses from 32 to 1024 nodes, which represents the upper bound for reasonable scalability given the modest problem size. Figure~\ref{fig:pelelmex_strong_scaling} shows the PeleLMeX strong scaling results for both the full timestep `Advance()' routine as well as the `advanceChemistry()' subroutine when using approach 1B. Here we see that the overall time step for the advanceChemistry routine ranges from 30 -- 15 seconds with a drop in scalability after 128 nodes (600,000 cells per node). The strong scaling results using approach 2B range from 15 seconds down to approximately 0.8 seconds. Both methods show an expected departure from the idealized scaling, but the improvement in speed when using 2B (modified Newton, analytic Jacobian) formulation is over an order of magnitude when compared with 1B (Newton-Krylov).

\subsection{Nyx}

\begin{figure}
\centering
\begin{tikzpicture}
\begin{loglogaxis}[
  title={Nyx Weak Scaling},
  xlabel={MPI Ranks},
  ylabel style={align=center},
  ylabel={Cells Evolved\\per Rank per Second},
  log basis x={10},
  legend columns=2,
  legend style={font=\tiny, at={(0.22,.75)},anchor=north west},
  legend cell align={left},
  width=0.95\linewidth,
  height=0.7\linewidth,
  every axis plot/.append style={thick},
  /pgf/number format/.cd,
        use comma,
        1000 sep={},
]

\addplot[thin, black!30, mark=none, dashed, forget plot] coordinates {
  (1,     8.07295544e+04)
  (13824, 8.07295544e+04)
};

\addplot[thin, black!30, mark=none, dashed, forget plot] coordinates {
  (1,     1.17159330e+06)
  (13824, 1.17159330e+06)
};

\addplot[thin, black!30, mark=none, dashed, forget plot] coordinates {
  (16,   1.55120658e+07)
  (2048, 1.55120658e+07)
};

\addplot[thin, black!30, mark=none, dashed, forget plot] coordinates {
  (16,    25885621.8382)
  (65536, 25885621.8382)
};

\addplot[blue, mark=triangle*] coordinates {
    (1,     8.07295544e+04)
    (216,   7.39997177e+04)
    (1728,  7.22532989e+04)
    (13824, 6.89399080e+04)
};

\addplot[red,, mark=square*] coordinates {
    (1,     1.17159330e+06)
    (216,   9.48401131e+05)
    (1728,  9.43600450e+05)
    (13824, 9.02000860e+05)
};

\addplot[black, mark=*] coordinates {
    (16,   1.55120658e+07)
    (32,   1.97789529e+07)
    (64,   1.83129288e+07) 
    (128,  1.56650520e+07) 
    (256,  1.89726096e+07)
    (512,  1.62941274e+07) 
    (1024, 1.35896211e+07) 
    (2048, 1.64947726e+07)
};

\addplot[magenta, mark=diamond*] coordinates {
    (16,    25885621.8382)
    (32,    21882619.3594431)
    (64,    22509190.0050277)
    (128,   24984866.9163928)
    (256,   22751336.1478055)
    (512,   21972749.2255819)
    (1024,  24789418.5101057)
    (2048,  21585613.3448272)
    (4096,  21768979.0801324)
    (8192,  23943964.2026122)
    (16384, 20501888.323182)
    (32768, 19776915.944495)
    (65536, 19724476.7556039)
};
\legend{Summit CPU, Summit GPU, Perlmutter, Frontier}
\end{loglogaxis}
\end{tikzpicture}
\caption{Weak scaling of the Lyman-$\alpha$ simulation with a random initial condition on Summit (1 to 2304 nodes) with CPUs (MPI + OMP) and GPUs (MPI + OMP + NVIDIA V100), Perlmutter (4 to 512 nodes, MPI + OMP + NVIDIA A100), and Frontier (2 to 8192 nodes, MPI + OMP + AMD MI250X). The gray dashed lines indicate the ideal scaling from the initial point.}
\label{fig:nyx_weak_scaling}
\end{figure}
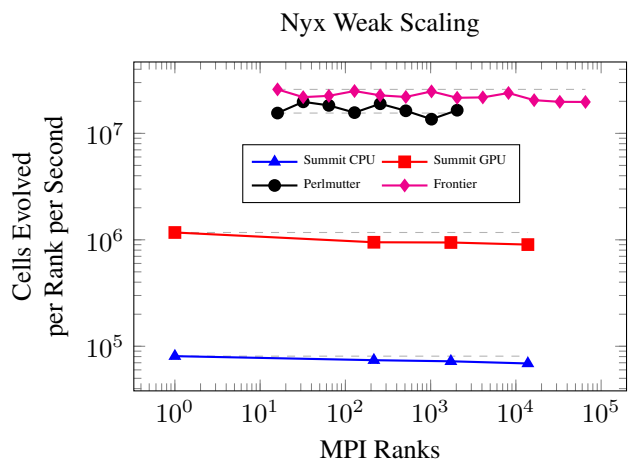

In Figure~\ref{fig:nyx_weak_scaling} we compare weak scaling results with Nyx on Summit, Perlmutter, and Frontier using the Lyman-$\alpha$ problem with an initial condition consisting of randomly placed particles while keeping the average density in each setup fixed.
For this problem, the metric of interest is the number of cells evolved per MPI rank per second (higher values indicate better performance).
Here, the number of cells is used to account for differing problem sizes which were chosen to be appropriate to the available memory. We can see that Frontier performance holds up well retaining 76\% of the throughput achieved on 2 nodes out to 8192 nodes where the problem size is $16384^3 \approx 4.39 \times 10^{12}$. This problem has about 4 boxes of $256^3$ cells per MPI rank on Frontier, with 1 rank per MI250X tile. At this problem size per node, Nyx has historically been unable to fit within the high bandwidth memory of 1 MPI rank per V100 on Summit. The results on Summit show similar levels of weak scaling efficiency while the Perlmutter retains at least 87\% of the throughput going from 4 nodes to 512 nodes.



\section{Conclusions}\label{s:conclusion}
In this article, we presented a batched strategy for integrating the independent ODEs that result from an operator split approach to solving coupled reacting flow PDEs in the context of the SUNDIALS time integration library, applications in cosmology and combustion, and modern GPU-based exascale supercomputers.
The PeleC and PeleLMeX combustion simulation codes and the Nyx cosmology simulation code each have some distinct features that provide different challenges (data ordering, tolerance selection, efficient linear solvers, GPU kernel launch overheads, low GPU utilization) and some that are common among all three codes (efficient memory usage, reduction operations).

In Pele applications, increased performance was achieved through several improvements: (i) distributing the block-structured AMR patches according to work estimates from the ODE solve, (ii) reordering data to allow coalescing of GPU memory accesses, (iii) selecting integrator tolerances dynamically based on typical values of the chemical species, and (iv) using batched linear solvers.
The improvements were tested together on the Frontier exascale supercomputer with the PeleLMeX code, and modest strong scaling was achieved up to 1024 nodes. The combustion simulation presented in Section~\ref{s:results}, demonstrates the use of the improved ODE integrator to simultaneously capture the broad range of spatial and temporal scales associated with the turbulent jet and their interaction with the detailed chemistry required to exhibit the low temperature $n$-dodecane ignition physics.
\addressed{MSD}{CSW}{Isn't there something stronger we can say here?  Was the faster chemistry enabling of the bigger runs?  Did the use of sundials contributed to a first ever run on an exascale system?}
Within the Nyx code, performance was improved by (i) increasing concurrency with OpenMP and GPU streams, (ii) fusing kernels within the CVODE integrator to increase the computational workload of kernels, and (iii) tiling loops over boxes in the block-structured AMR to better use the high-bandwidth GPU memory.
These improvements were tested on the Summit, Perlmutter, and Frontier supercomputers and showed excellent weak scaling up to 8192 nodes.
The Pele and Nyx applications both benefited from the use of memory pools that were accessible by the application code, AMReX framework, and SUNDIALS to prevent resource contention and reduce data movement.
The codes also benefited from faster GPU reduction operations within SUNDIALS. Scientific outcomes related to these code advances are varied. They were used in Lyman-$\alpha$ forest comparison studies at meaningful simulation resolutions which considered the effect on the accuracy of the different temperature profiles on the overall simulation results \citep{Chabanier2023}. They includes the largest hydrodynamical simulations to date of the Lyman-$\alpha$ forest with high physical resolution \citep{Chabanier:inprep}.
\addressed{CSW}{}{Isn't there something stronger we can say here?  Was the faster ODE solver enabling of the bigger runs?  Did the use of sundials contributed to a first ever run on an exascale system? Cody added a sentence to the beginning of the next paragraph along these lines. Maybe there is still something more concrete we can say.}

The results not only show that batching of independent ODEs can be an effective strategy for operator split multiphysics applications, especially those with reacting flows, but also that the Pele and Nyx applications and the SUNDIALS time integration library are ready for the exascale era of high performance computing.
Furthermore it is noteworthy that the collaboration of the SUNDIALS, Pele, and Nyx teams enabled by the Exascale Computing Project was effective.
It can serve as encouragement for further tight collaborations of those working on math libraries and science application teams in the future.

\begin{acks}
This research used resources of the Oak Ridge Leadership Computing Facility, which is a DOE Office of Science User Facility supported under Contract DE-AC05-00OR22725.
This research also used resources of the National Energy Research Scientific Computing Center (NERSC), a U.S. Department of Energy Office of Science User Facility located at Lawrence Berkeley National Laboratory, operated under Contract No. DE-AC02-05CH11231 using NERSC awards ERCAP0028634 and ERCAP0026830.
A portion of the research was performed using computational resources sponsored by the Department of Energy’s Office of Energy Efficiency and Renewable Energy and located at the National Renewable Energy Laboratory.
This work was performed under the auspices of the U.S. Department of Energy by Lawrence Livermore National Laboratory under Contract DE-AC52-07NA27344. Lawrence Livermore National Security, LLC. LLNL-JRNL-860348.
This work was authored in part by the National Renewable Energy Laboratory, operated by Alliance for Sustainable Energy, LLC, for the U.S. Department of Energy (DOE) under Contract No. DE-AC36-08GO28308.
Research at LBNL was provided by the U.S.~Department of Energy, Office of Science, Office of Advanced Scientific Computing Research, Exascale Computing Project under contract DE-AC02-05CH11231.
The views expressed in the article do not necessarily represent the views of the DOE or the U.S. Government. The U.S. Government retains and the publisher, by accepting the article for publication, acknowledges that the U.S. Government retains a nonexclusive, paid-up, irrevocable, worldwide license to publish or reproduce the published form of this work, or allow others to do so, for U.S. Government purposes.

\end{acks}

\begin{dci}
The authors declare no conflicts of interests.
\end{dci}

\begin{funding}
This research was supported by the Exascale Computing Project (17-SC-20-SC), a collaborative effort of the U.S. Department of Energy Office of Science and the National Nuclear Security Administration. Support for this work was also provided in part by the U.S. Department of Energy, Office of Science, Office of Advanced Scientific Computing Research, Scientific Discovery through Advanced Computing (SciDAC) Program through the Frameworks, Algorithms, and Scalable Technologies for Mathematics (FASTMath) Institute, under Lawrence Livermore National Laboratory subcontract B626484 and DOE award DE-SC0021354.
\end{funding}

\bibliographystyle{SageH}
\bibliography{references}

\end{document}